\documentclass[12pt]{article}
\usepackage{latexsym, amssymb, amsmath, amscd, amsfonts, epsfig, graphicx, colordvi,verbatim,ifpdf,extarrows}
\usepackage{amsfonts, amsmath, amssymb}
\usepackage{amssymb,amsfonts,amsmath,latexsym,epsfig,cite, psfrag,eepic,color}
\usepackage{amscd,graphics}
\usepackage{latexsym, amssymb,  amsmath,amscd, amsfonts, epsfig, graphicx, colordvi,amsthm}

\usepackage{graphicx}
\usepackage{epstopdf}
\usepackage{color}
\usepackage{ifpdf}
\usepackage{fancybox}
\usepackage[font=small,labelfont=bf,labelsep=none]{caption}
\usepackage{float}

\newtheorem{thm}{Theorem}[section]
\newtheorem{prop}[thm]{Proposition}
\newtheorem{conj}[thm]{Conjecture}

\newtheorem{defi}[thm]{Definition}
\newtheorem{lem}[thm]{Lemma}

\def\pf{\noindent{\it Proof.} }
\def\qed{\nopagebreak\hfill{\rule{4pt}{7pt}}
\medbreak}

\numberwithin{equation}{section}

\def\qed{\nopagebreak\hfill{\rule{4pt}{7pt}}
\medbreak}

\newlength{\boxedparwidth}
  {\begin{center} \begin{tabular}{|@{\hspace{.315in}}c@{\hspace{.15in}}|}
                  \hline \\ \begin{minipage}[t]{\boxedparwidth}
                  \setlength{\parindent}{.25in}}%
  {\end{minipage} \\ \\ \hline \end{tabular} \end{center}}

\begin{document}

\begin{center}

 {\Large \bf An overpartition analogue of Bressoud conjecture for even moduli}
\end{center}

\begin{center}
 {Y.H. Chen}$^{1}$, {T.T. Gu}$^{2}$, {Thomas Y. He}$^{3}$, {F. Tang}$^{4}$ and
  {J.J. Wei}$^{5}$ \vskip 2mm

$^{1,2,3,4,5}$ School of Mathematical Sciences, Sichuan Normal University, Chengdu 610066, P.R. China

   \vskip 2mm

  $^1$chenyh@stu.sicnu.edu.cn, $^2$guyichen0816@163.com, $^3$heyao@sicnu.edu.cn,  $^4$tangfan@stu.sicnu.edu.cn,  $^5$wei@stu.sicnu.edu.cn
\end{center}

\noindent{\bf Abstract.}  In  1980,  Bressoud conjectured a combinatorial identity $A_j=B_j$ for $j=0$ or $1$.
In this paper,  we introduce  a  new partition function $\overline{B}_0$ which can be viewed as an overpartition analogue of the partition function $B_0$. An overpartition is a
partition such that the last occurrence of a part can be overlined.
We build a bijection to get a relationship between $\overline{B}_0$ and $B_1$, based on which an overpartition analogue of Bressoud's conjecture for $j=0$ is obtained.



\section{Introduction}

A partition $\pi$ of a positive integer $n$ is a finite non-increasing sequence of positive integers $\pi=(\pi_1,\pi_2,\ldots,\pi_\ell)$ such that $\sum_{i=1}^{\ell}\pi_i=n.$ The $\pi_i$ are called the parts of $\pi$.  Let  $\ell(\pi)$ be the number of parts of $\pi$ and  $|\pi|$ be the sum of parts of $\pi$.

 Throughout this paper, we assume that $\alpha_1,\alpha_2,\ldots, \alpha_\lambda$ and $\eta$ are integers such that
\begin{equation*}\label{cond-alpha}
0<\alpha_1<\cdots<\alpha_\lambda<\eta, \quad \text{and} \quad \alpha_i=\eta-\alpha_{\lambda+1-i}\quad \text{for} \quad 1\leq i\leq \lambda.
\end{equation*}

In 1980, Bressoud  \cite{Bressoud-1980} introduced the following two partition functions.

\begin{defi}[Bressoud] For $j=0$ or $1$  and $(2k+j)/2> r\geq\lambda\geq0$,   define the partition function $A_j(\alpha_1,\ldots,\alpha_\lambda;\eta,k,r;n)$ to be the number of partitions of $n$ into parts congruent to $0,\alpha_1,\ldots,\alpha_\lambda\pmod\eta$ such that
\begin{itemize}
\item[\rm{(1)}] If $\lambda$ is even, then only multiples of $\eta$ may be repeated and no part is congruent to $0,\pm\eta(r-\lambda/2)  \pmod{\eta(2k-\lambda+j)}${\rm{;}}

    \item[\rm{(2)}] If $\lambda$ is odd and $j=1$, then only multiples of ${\eta}/{2}$ may be repeated, no part is congruent to $\eta\pmod{2\eta}$, and no part is congruent to $0,\pm{\eta}(2r-\lambda)/{2} \pmod {\eta(2k-\lambda+1)}${\rm{;}}

        \item[\rm{(3)}] If $\lambda$ is odd and $j=0$, then only multiples of ${\eta}/{2}$ which are not congruent to ${\eta}(2k-\lambda)/{2}\pmod{\eta(2k-\lambda)}$ may be repeated, no part is congruent to $\eta\pmod{2\eta}$, no part is congruent to $0\pmod{2\eta(2k-\lambda)}$, and no part is congruent to $\pm{\eta}(2r-\lambda)/{2} \pmod {\eta(2k-\lambda)}$.
  \end{itemize}
  \end{defi}

\begin{defi}[Bressoud] \label{Bress-B-function} For $j=0$ or $1$  and $k\geq r\geq\lambda\geq0$,  define the partition function
$B_j(\alpha_1,\ldots,\alpha_\lambda;\eta,k,r;n)$  to be the number of partitions  $\pi=(\pi_1,\pi_2,\ldots,\pi_\ell)$  of $n$  satisfying the following conditions{\rm{:}}
\begin{itemize}
\item[{\rm (1)}] For $1\leq i\leq \ell$, $\pi_i\equiv0,\alpha_1,\ldots,\alpha_\lambda\pmod{\eta}${\rm{;}}

\item[{\rm (2)}] Only multiples of $\eta$ may be repeated{\rm{;}}

\item[{\rm (3)}] For $1\leq i\leq \ell-k+1$, $ \pi_i\geq\pi_{i+k-1}+\eta$  with strict inequality if $\eta\mid\pi_i${\rm{;}}

\item[{\rm (4)}] At most $r-1$ of the $\pi_i$ are less than or equal to $\eta${\rm{;}}

\item[{\rm (5)}] For $1\leq i\leq \ell-k+2$, if $\pi_i\leq\pi_{i+k-2}+\eta$ with strict inequality if $\eta\nmid\pi_i$,  then \[[\pi_i/\eta]+\cdots+[\pi_{i+k-2}/\eta]\equiv r-1+V_\pi(\pi_i)\pmod{2-j},\]
  where $V_\pi(N)$ denotes the number of parts not exceeding $N$ which are not divisible by $\eta$ in $\pi$ and $[\ ]$ denotes the greatest integer function.
  \end{itemize}
 \end{defi}

 Bressoud  \cite{Bressoud-1980} made the following conjecture.

  \begin{conj}[Bressoud] \label{Bressoud-conjecture-j} For $j=0$ or $1$, $(2k+j)/2> r\geq\lambda\geq0$ and $n\geq 0$,
 \begin{equation*}\label{Bressoud-conj-1}
A_j(\alpha_1,\ldots,\alpha_\lambda;\eta,k,r;n)=B_j(\alpha_1,\ldots,\alpha_\lambda;\eta,k,r;n).
 \end{equation*}
\end{conj}

This conjecture specializes to a wide variety of well-known theorems, including the Rogers-Ramanujan-Gordon theorem \cite{Gordon-1961}, the Bressoud-Rogers-Ramanujan theorem \cite{Bressoud-1979} and the Andrews-G\"ollnitz-Gordon theorem \cite{Andrews-1967}. Bressoud \cite{Bressoud-1980} proved Conjecture \ref{Bressoud-conjecture-j} for $\lambda=1$ and $j=0$, which has been called the Bressoud-G\"ollnitz-Gordon theorem. Conjecture \ref{Bressoud-conjecture-j}   for $\eta=\lambda+1$ and $j=1$ was verified by Andrews \cite{Andrews-1974m}. Kim and Yee \cite{Kim-Yee-2014} showed that Conjecture \ref{Bressoud-conjecture-j} holds for $j=1$ and $\lambda=2$.

Recently, Kim  \cite{Kim-2018} proved Bressoud's conjecture for $j=1$. Bressoud's conjecture for $j=0$  was resolved by He, Ji and Zhao \cite{He-Ji-Zhao}. To this end, Kim \cite{Kim-2018} and He, Ji and Zhao \cite{He-Ji-Zhao} established the following theorem for $j=1$ and $j=0$ respectively.

\begin{thm}\label{berssoud-kim-1}
For $j=0$ or $1$ and $(2k+j)/2> r\geq\lambda\geq0$,
\begin{equation*}\label{Bressoud-conj-defi-e-11}
 \begin{split}
 &\sum_{n\geq0}B_j(\alpha_1,\ldots,\alpha_\lambda;\eta,k,r;n)q^n\\[5pt]
 &=\frac{(-q^{\alpha_1},\ldots,-q^{\alpha_\lambda};q^{\eta})_\infty
 (q^{\eta(r-\frac{\lambda}{2})},q^{\eta(2k-r-\frac{\lambda}{2}+j)},
 q^{\eta(2k-\lambda+j)};q^{\eta(2k-\lambda+j)})_\infty}{(q^\eta;q^\eta)_\infty}.
 \end{split}
 \end{equation*}
 \end{thm}

 From now on, we assume that $|q|<1$ and adopt the standard notation \cite{Andrews-1976}:
\[(a;q)_\infty=\prod_{i=0}^{\infty}(1-aq^i), \quad (a;q)_n=\frac{(a;q)_\infty}{(aq^n;q)_\infty},\]
and
\[(a_1,a_2,\ldots,a_m;q)_\infty=(a_1;q)_\infty(a_2;q)_\infty\cdots(a_m;q)_\infty.
\]

Overpartitions were first introduced by  \cite{Corteel-Lovejoy-2004}. In recent years, overpartitions have been heavily studied. He, Ji and Zhao gave overpartitions analogues of Bressoud's conjecture for $j=1$ and $j=0$ in \cite{he-ji-zhao-2022} and \cite{He-Ji-Zhao} respectively. In \cite{he-ji-zhao-2022,He-Ji-Zhao}, an overpartition is a
partition such that the first occurrence of a part can be overlined and the parts in an overpartition are ordered  as follows: $1<\bar{1}<2<\bar{2}<\cdots.$

In this paper, an overpartition is defined to be a
partition such that the last occurrence of a part can be overlined. We impose the following order on the parts of an overpartition:
\begin{equation}\label{order}
\bar{1}<{1}<\bar{2}<{2}<\cdots.
\end{equation}
Throughout this paper,
we adopt the following convention. For positive integer $t$, we define $|t|$ (resp. $|\bar{t}|$) to be the size of $t$  (resp. $\bar{t}$), that is, $|t|=t$ (resp. $|\bar{t}|=t$). For positive integer $b$, we define $t\pm b$ (resp. $\overline{t}\pm b$) as a non-overlined part (resp. an overlined part) of size  $t\pm b$. Let $\pi$ be an overpartition. We use $f_t(\pi)$ (resp. $f_{\bar{t}}(\pi)$) to denote the number of parts equal to $t$ (resp. $\bar{t}$) in $\pi$. We use $f_{\leq \eta}(\pi)$ to denote the number of parts not exceeding $\eta$ in $\pi$.

The main objective of this paper is to give a new overpartition analogue of Bressoud's conjecture for $j=0$,
 involving the following two partition functions $\overline{A}_0(\alpha_1,\ldots,\alpha_\lambda;\eta, k,r;n)$ and $\overline{B}_0(\alpha_1,\ldots,\alpha_\lambda;\eta, k,r;n)$. The former partition function was introduced by He, Ji and Zhao \cite[Definition 1.15]{he-ji-zhao-2022}.

 \begin{defi}[He-Ji-Zhao]\label{defi-O-A}
 For $k\geq r\geq \lambda\geq0$, define  $\overline{A}_0(\alpha_1,\ldots,\alpha_\lambda;\eta,k,r;n)$ to be the number of overpartitions of $n$ satisfying
$\pi_i\equiv0,\alpha_1,\ldots,\alpha_\lambda\pmod{\eta}$ such that
 \begin{itemize}
 \item[{\rm (1)}] If $\lambda$ is even, then only multiplies of $\eta$ may be non-overlined and  there is no non-overlined part congruent to $0,\pm\eta(r-\lambda/2) \pmod {\eta(2k-\lambda-1)}${\rm;}

 \item[{\rm (2)}]  If $\lambda$ is odd, then only multiples of  ${\eta}/{2}$  may be non-overlined, no non-overlined part is congruent to $\eta \pmod{2\eta}$,  no non-overlined part is congruent to $0,\pm{\eta}(2r-\lambda)/{2} \pmod {\eta(2k-\lambda-1)}$, and no overlined part is congruent to ${\eta}/{2}\pmod \eta$.
  \end{itemize}
    \end{defi}

     \begin{defi}\label{defi-O-B} For $k>r\geq \lambda\geq0$, define $\overline{B}_0(\alpha_1,\ldots,\alpha_\lambda;\eta,k,r;n)$ to be the number of overpartitions $\pi=(\pi_1,\pi_2,\ldots,\pi_\ell)$ of $n$   subject to the following conditions{\rm{:}}
 \begin{itemize}
  \item[{\rm (1)}] For $1\leq i\leq \ell$, $\pi_i\equiv0,\alpha_1,\ldots,\alpha_\lambda\pmod{\eta}${\rm{;}}

 \item[{\rm (2)}] Only  multiples of $\eta$ may be non-overlined{\rm{;}}

 \item[{\rm (3)}]  For $1\leq i\leq \ell-k+1$, $\pi_i\geq\pi_{i+k-1}+\eta$ with strict inequality if  $\pi_i$  is non-overlined{\rm{;}}

 \item[{\rm (4)}] $f_{\leq \eta}(\pi)\leq r${\rm{;}}

 \item[{\rm (5)}] If $f_{\leq \eta}(\pi)=r$ and $\overline{\eta}$ does not occur in $\pi$, then there exists $i$ such that $1\leq i\leq \ell-k+2$, $\pi_i<\overline{2\eta}$, and $\pi_i\leq\pi_{i+k-2}+\eta$ with strict inequality if   $\pi_i$ is overlined{\rm{;}}

 \item[{\rm (6)}] For $1\leq i\leq \ell-k+2$, if $\pi_i\leq\pi_{i+k-2}+\eta$ with strict inequality if   $\pi_i$ is overlined, then
 \[\left[|\pi_i|/\eta\right]+\cdots+\left[|\pi_{i+k-2}|/\eta\right]\equiv r-1+\overline{V}_{\pi}(\pi_i)+\overline{O}_{\pi}(\pi_{i+k-2})\pmod{2},\]
 where $\overline{V}_{\pi}(\pi_i)$ denotes the number of overlined parts not exceeding $\pi_i$ which are not divisible by $\eta$ in $\pi$ and $\overline{O}_{\pi}(\pi_{i+k-2})$ denotes the number of overlined parts greater than or equal to $\pi_{i+k-2}$ which are divisible by $\eta$ in $\pi$.
\end{itemize}
\end{defi}

Let $\pi=(\pi_1,\pi_2,\ldots,\pi_\ell)$ be an overpartition satisfying the conditions (1)-(3) in Definition \ref{defi-O-B}. Assume that there exists  $i$ such that $1\leq i\leq \ell-k+2$, $\pi_i<\overline{2\eta}$, and $\pi_i\leq\pi_{i+k-2}+\eta$ with strict inequality if   $\pi_i$ is overlined. It follows from the condition (3) in Definition \ref{defi-O-B} that $\pi_{i+k-1}\leq \pi_i-\eta<\overline{2\eta}-\eta=\overline{\eta}$ and $\pi_{i-1}\geq \pi_{i+k-2}+\eta>\eta$. So, we get
\[
\begin{split}
&\quad\left[|\pi_i|/\eta\right]+\cdots+\left[|\pi_{i+k-2}|/\eta\right]\\
&=\left[|\pi_i|/\eta\right]+\cdots+\left[|\pi_{i+k-2}|/\eta\right]+\left[|\pi_{i+k-1}|/\eta\right]+\cdots+\left[|\pi_{\ell}|/\eta\right]\\
&=\overline{V}_{\pi}(\pi_i)-f_{\leq \eta}(\pi)+2(f_{\overline{\eta}}(\pi)+f_{\eta}(\pi))\\
&\equiv f_{\leq \eta}(\pi)+\overline{V}_{\pi}(\pi_i) \pmod 2.
\end{split}
\]

Hence, for an overpartition $\pi$ counted by $\overline{B}_0(\alpha_1,\ldots,\alpha_\lambda;\eta,k,r;n)$ without overlined parts divisible by $\eta$, we have $f_{\leq \eta}(\pi)\neq r$, and so $f_{\leq \eta}(\pi)\leq r-1$. Furthermore, if we change the overlined parts in $\pi$ to non-overlined parts, then we get an ordinary partition counted by ${B}_0(\alpha_1,\ldots,\alpha_\lambda;\eta,k,r;n)$. Thus, we  say that $\overline{B}_0(\alpha_1,\ldots,\alpha_\lambda;\eta,k,r;n)$ can be considered as an overpartition analogue of  ${B}_0(\alpha_1, \ldots,\alpha_\lambda;\eta,k,r;n)$.

The following theorem is the main result of this paper, which can be regarded as a new overpartition analogue of Bressoud's conjecture for $j=0$.
\begin{thm}\label{main-Bressoud-conjecture-0} For
 $k> r\geq\lambda\geq0,$ $k-1>\lambda$  and $n\geq 0,$
 \begin{equation*}\label{eqn-Bressoud-conj-0}
\overline{A}_0(\alpha_1,\ldots,\alpha_\lambda;\eta,k,r;n)=\overline{B}_0(\alpha_1,\ldots,\alpha_\lambda;\eta,k,r;n).
 \end{equation*}
\end{thm}

The inspiration of the proof of Theorem \ref{main-Bressoud-conjecture-0} comes from \cite[Section 5]{he-ji-zhao-2022}. In Section 2, we give an equivalent statement of Theorem \ref{main-Bressoud-conjecture-0}, which is stated in Theorem \ref{eqv-main-0}. Section 3 is a
preliminary for the proof of Theorem \ref{eqv-main-0}, in which we introduce some necessary notions and results about $(k-1)$-bands.   In Section 4, we define the $(k-1)$-reduction operation and the $(k-1)$-augmentation operation, which allow us to provide a combinatorial proof of Theorem \ref{eqv-main-0}.

\section{Equivalent statement of Theorem \ref{main-Bressoud-conjecture-0}}

 The generating function of $\overline{A}_0(\alpha_1,\ldots,\alpha_\lambda;\eta,k,r;n)$ is given by He, Ji and Zhao \cite[Theorem 1.18]{he-ji-zhao-2022}.
 \begin{thm}[He-Ji-Zhao]\label{eqn-bressoud-over-0}
 For $k\geq r\geq \lambda,$
\begin{equation*}\label{gf-overlineA_0-eq}
 \begin{split}
 & \sum_{n\geq0}\overline{A}_{0}(\alpha_1,\ldots,\alpha_\lambda;\eta,k,r;n)q^n\\[5pt]
  &=\frac{(-q^{\alpha_1},\ldots,-q^{\alpha_\lambda},-q^{\eta};q^{\eta})_\infty(q^{\eta(r-\frac{\lambda}{2})},q^{\eta
 (2k-r-\frac{\lambda}{2}-1)},q^{\eta(2k-\lambda-1)}
 ;q^{\eta(2k-\lambda-1)})_\infty}{(q^\eta;q^\eta)_\infty}.
 \end{split}
 \end{equation*}
 \end{thm}

 It follows from Theorems \ref{berssoud-kim-1} and \ref{eqn-bressoud-over-0} that Theorem \ref{main-Bressoud-conjecture-0} is equivalent to the following relationship between $\overline{B}_0$ and $B_1$.
 \begin{thm}\label{b-0-over-000-thm}
 For $k> r\geq\lambda\geq0$ and $k-1>\lambda,$
 \begin{equation*}\label{b-0-over-000}
 \begin{split}
 &\sum_{n\geq0}\overline{B}_0(\alpha_1,\ldots,\alpha_\lambda;
 \eta,k,r;n)q^n=(-q^\eta;q^\eta)_\infty
 \sum_{n\geq0}B_1(\alpha_1,\ldots,\alpha_\lambda;
 \eta,k-1,r;n)q^n.
 \end{split}
\end{equation*}
\end{thm}

To give a combinatorial statement of Theorem \ref{b-0-over-000-thm}, we prepare the following sets. Here and in the remaining of this paper, we assume that $k$, $r$ and $\lambda$ are integers such that
$k> r\geq\lambda\geq0$ and $k-1>\lambda$.
\begin{itemize}

\item $\overline{\mathcal{B}}_0(\alpha_1,\ldots,\alpha_\lambda;\eta,k,r)$:  the set of overpartitions  counted by $\overline{B}_0(\alpha_1,\ldots,\alpha_\lambda;\eta,k,r;n)$ for  $n\geq 0$;

\item $\overline{\mathcal{B}}(\alpha_1,\ldots,\alpha_\lambda;\eta,k,r)$: the set of overpartitions satisfying the conditions (1)-(4) in Definition \ref{defi-O-B};

\item ${\mathcal{B}}_1(\alpha_1,\ldots,\alpha_\lambda;\eta,k,r)$: the set of overpartitions $\pi$ in  $\overline{\mathcal{B}}(\alpha_1,\ldots,\alpha_\lambda;\eta,k,r)$ such that there are no overlined parts divisible by $\eta$ in $\pi$ and $f_{\leq \eta}(\pi)\leq r-1$;

\item $\mathcal{D}_\eta$:  the set of  partitions with distinct parts divisible by $\eta$.

\end{itemize}

Let $\pi$ be an overpartition in ${\mathcal{B}}_1(\alpha_1,\ldots,\alpha_\lambda;\eta,k,r)$ with $|\pi|=n$. If we change the overlined parts in $\pi$ to non-overlined parts, then we get an ordinary partition counted by  $B_1(\alpha_1,\ldots,\alpha_\lambda;\eta,k,r;n)$, and vice versa. Hence, the combinatorial statement of Theorem \ref{b-0-over-000-thm} can be stated as follows, which is an equivalent statement of Theorem \ref{main-Bressoud-conjecture-0}.

\begin{thm}\label{eqv-main-0}
There is a bijection $\Phi$ between  $\mathcal{\overline{B}}_0(\alpha_1,\ldots,\alpha_\lambda;\eta,k,r)$ and $\mathcal{D}_\eta\times{\mathcal{B}}_{1}(\alpha_1,\ldots,\alpha_\lambda;$\break$ \eta,k-1,r)$, namely, for an overpartition $\pi \in\mathcal{\overline{B}}_0(\alpha_1,\ldots,\alpha_\lambda;\eta,k,r)$, we have $\Phi(\pi)=(\tau,\mu) \in \mathcal{D}_\eta\times\mathcal{{B}}_{1}(\alpha_1,\ldots,\alpha_\lambda;\eta,k-1,r) $ such that  $|\pi|=|\tau|+|\mu|$ and $\ell(\pi)=\ell(\tau)+\ell(\mu)$.
\end{thm}

\section{$(k-1)$-bands}

In this section, we introduce the definition of $m$-bands and investigate properties of $m$-bands. We will mainly focus on $(k-1)$-bands. It is worth mentioning that the definition of $(k-1)$-bands is almost the same as that in \cite[Section 3]{he-ji-zhao-2022} except that the parts in an overpartition are ordered as in \eqref{order} in this paper.

Let $\pi$ be an overpartition in $\overline{\mathcal{B}}(\alpha_1,\ldots,\alpha_\lambda;\eta,k,r)$.
For $m\geq 1$, we use $\{\pi_{i+l}\}_{0\leq l\leq m-1}$ to denote the set of the $m$ consecutive parts  $\pi_i,\pi_{i+1},\ldots,\pi_{i+m-1}$ of $\pi$. We say that $\{\pi_{i+l}\}_{0\leq l\leq m-1}$ is a $m$-band of $\pi$ if $\pi_i\leq\pi_{i+m-1}+\eta$ with strict inequality if   $\pi_i$ is overlined, in which we have $m\leq k-1$.

For $t\geq 1$, we say that a $m$-band $\{\pi_{i+l}\}_{0\leq l\leq m-1}$ belongs to $[(t-1)\eta,(t+1)\eta]$ (resp. $[(t-1)\eta,\overline{(t+1)\eta})$) if $\pi_{i+m-1}\geq (t-1)\eta$ and $\pi_i\leq (t+1)\eta$ (resp. $\pi_i<\overline{(t+1)\eta}$). Then, we have
\begin{lem}\label{t+1}
For $t\geq 1$, let $\pi$ and $\mu$ be two overpartitions in $\overline{\mathcal{B}}(\alpha_1,\ldots,\alpha_\lambda;\eta,k,r)$ such that $\pi$ can be obtained by inserting ${t\eta}$ as a non-overlined part into $\mu$. Then,
\begin{itemize}
\item[{\rm(1)}] there are no $(k-1)$-bands of $\mu$ belonging to $[(t-1)\eta,(t+1)\eta]$;

\item[{\rm(2)}] $\{\pi_{i+l}\}_{0\leq l\leq k-2}$ is a $(k-1)$-band of $\pi$ belonging to $[(t-1)\eta,(t+1)\eta]$ {\rm(}resp. $[(t-1)\eta,\overline{(t+1)\eta})${\rm)} if and only if $\{\mu_{i+l}\}_{0\leq l\leq k-3}$ is a $(k-2)$-band of $\mu$ belonging to $[(t-1)\eta,(t+1)\eta]$ {\rm(}resp. $[(t-1)\eta,\overline{(t+1)\eta})${\rm)}. Moreover, we have
    \begin{equation}\label{cong-m-band}
    \begin{split}
&\quad [|\pi_{i}|/\eta]+\cdots+[|\pi_{i+k-2}|/\eta]+\overline{V}_{\pi}(\pi_{i})+\overline{O}_{\pi}(\pi_{i+k-2})\\
&\equiv [|\mu_{i}|/\eta]+\cdots+[|\mu_{i+k-3}|/\eta]+\overline{V}_{\mu}(\mu_{i})+\overline{O}_{\mu}(\mu_{i+k-3})+t\pmod2.
\end{split}
\end{equation}

\end{itemize}
\end{lem}

\pf (1). Suppose to the contrary that there exists a $(k-1)$-band $\{\mu_{s+l}\}_{0\leq l\leq k-2}$ of $\mu$ belonging to $[(t-1)\eta,(t+1)\eta]$. Then, we have  $\mu_{s+k-2}\geq (t-1)\eta$ and $\mu_s\leq (t+1)\eta$.
Using the condition (3) in Definition \ref{defi-O-B}, we see that 
\[\mu_{s-1}\geq \mu_{s+k-2}+\eta\text{ with strict inequality if }\mu_{s-1}\text{ is non-overlined,}\]
and 
\[\mu_{s+k-1}\leq \mu_{s}-\eta\text{ with strict inequality if }\mu_{s}\text{ is non-overlined.}\]
It yields that  $\mu_{s-1}>t\eta$ and $\mu_{s+k-1}<t\eta$.
Therefore, the consecutive $k$ parts $\pi_s,\ldots,\pi_{s+k-1}$ of $\pi$ consist of the parts $\mu_s,\ldots,\mu_{s+k-2}$ of $\mu$ together with the new inserted part ${t\eta}$. It can be checked that  $\{\pi_{s+l}\}_{0\leq l\leq k-1}$ is a $k$-band of $\pi$, which leads to a contradiction. Hence, there are no $(k-1)$-bands of $\mu$ belonging to $[(t-1)\eta,(t+1)\eta]$.

(2). Assume that $\{\mu_{i+l}\}_{0\leq l\leq k-3}$ is a $(k-2)$-band of $\mu$ belonging to $[(t-1)\eta,(t+1)\eta]$ (resp. $[(t-1)\eta,\overline{(t+1)\eta})$). By the condition (1), we know that there are no $(k-1)$-bands of $\mu$ belonging to $[(t-1)\eta,(t+1)\eta]$. With the similar argument in the proof of the condition (1), we obtain that the consecutive $k-1$ parts $\pi_i,\ldots,\pi_{i+k-2}$ of $\pi$ consist of the parts $\mu_i,\ldots,\mu_{i+k-3}$ of $\mu$ together with the new inserted part ${t\eta}$, and so $\{\pi_{i+l}\}_{0\leq l\leq k-2}$ is a $(k-1)$-band of $\pi$ belonging to $[(t-1)\eta,(t+1)\eta]$ (resp. $[(t-1)\eta,\overline{(t+1)\eta})$). This completes the proof of the sufficiency.

Conversely, assume that $\{\pi_{i+l}\}_{0\leq l\leq k-2}$ is a $(k-1)$-band of $\pi$ belonging to $[(t-1)\eta,(t+1)\eta]$ (resp. $[(t-1)\eta,\overline{(t+1)\eta})$). Assume that the new inserted part ${t\eta}$ is the $g$-th part of $\pi$, that is, $\pi_g=t\eta$.

We proceed to show that $i\leq g\leq i+k-2$. Suppose to the contrary that $g<i$ or $g>i+k-2$. In view of the condition (3) in Definition \ref{defi-O-B}, we get $\pi_{i+k-2}<\pi_g-\eta=(t-1)\eta$ or $\pi_i>\pi_g+\eta=(t+1)\eta$, which leads to a contradiction. Hence, we have $i\leq g\leq i+k-2$.

Recall that $\pi$ can be obtained by inserting ${t\eta}$ as a non-overlined part into $\mu$, so we have $\mu_{i+l}=\pi_{i+l}$ for $0\leq l<g-i$ and $\mu_{i+l}=\pi_{i+l+1}$ for $g-i\leq l\leq k-3$. So, $\{\mu_{i+l}\}_{0\leq l\leq k-3}$ is a $(k-2)$-band of $\mu$ belonging to $[(t-1)\eta,(t+1)\eta]$ (resp. $[(t-1)\eta,\overline{(t+1)\eta})$).

Again by the construction of $\pi$, we obtain that  $\overline{V}_{\pi}(\pi_{i})=\overline{V}_{\mu}(\mu_{i})$,  $\overline{O}_{\pi}(\pi_{i+k-2})=\overline{O}_{\mu}(\mu_{i+k-3})$, and
\[[|\pi_{i}|/\eta]+\cdots+[|\pi_{i+k-2}|/\eta]=[|\mu_{i}|/\eta]+\cdots+[|\mu_{i+k-3}|/\eta]+t,\]
which leads to \eqref{cong-m-band}. This completes the proof.
\qed

With the similar argument in the proof of Lemma \ref{t+1}, we can get
\begin{lem}\label{t+1-o}
For $t\geq 1$, let $\pi$ and $\mu$ be two overpartitions in $\overline{\mathcal{B}}(\alpha_1,\ldots,\alpha_\lambda;\eta,k,r)$ such that $\pi$ can be obtained by inserting $\overline{t\eta}$ as an overlined part into $\mu$. Then,
\begin{itemize}
\item[{\rm(1)}] there are no $(k-1)$-bands of $\mu$ belonging to $[(t-1)\eta,\overline{(t+1)\eta})$;

\item[{\rm(2)}] $\{\pi_{i+l}\}_{0\leq l\leq k-2}$ is a $(k-1)$-band of $\pi$ belonging to $[(t-1)\eta,\overline{(t+1)\eta})$ if and only if $\{\mu_{i+l}\}_{0\leq l\leq k-3}$ is a $(k-2)$-band of $\mu$ belonging to $[(t-1)\eta,\overline{(t+1)\eta})$. Moreover, we have
    \begin{equation}\label{cong-m-band-o}
    \begin{split}
&\quad [|\pi_{i}|/\eta]+\cdots+[|\pi_{i+k-2}|/\eta]+\overline{V}_{\pi}(\pi_{i})+\overline{O}_{\pi}(\pi_{i+k-2})\\
&\equiv[|\mu_{i}|/\eta]+\cdots+[|\mu_{i+k-3}|/\eta]+\overline{V}_{\mu}(\mu_{i})+\overline{O}_{\mu}(\mu_{i+k-3})+t+1\pmod2.
\end{split}
\end{equation}

\end{itemize}
\end{lem}

For a $(k-1)$-band $\{\pi_{i+l}\}_{0\leq l\leq k-2}$ of $\pi$, if $\{\pi_{i+l}\}_{0\leq l\leq k-2}$ satisfies the congruence condition
\[\left[|\pi_i|/\eta\right]+\cdots+\left[|\pi_{i+k-2}|/\eta\right]\equiv r-1+\overline{V}_{\pi}(\pi_i)+\overline{O}_{\pi}(\pi_{i+k-2})\pmod{2},\]
then we say that $\{\pi_{i+l}\}_{0\leq l\leq k-2}$ is even. Otherwise, we say that it is odd.

For example, let $\pi$ be an overpartition in $\overline{\mathcal{B}}(3,5,7;10,5,4)$  given by
\begin{equation*}\label{mark-exa-1}
\begin{split}
&\pi=(80,80,\overline{80},70,\overline{70},\overline{67},
60,\overline{60},\overline{55},\overline{53},\overline{50},\overline{47},\overline{45},\overline{43},\overline{37},\overline{35},\\[3pt]
&\ \ \ \ \ \ \ \overline{27},{20},20,\overline{20},\overline{13},\overline{10},\overline{7},
\overline{5},\overline{3}).
\end{split}
\end{equation*}
There are ten $4$-bands of $\pi$.
  \[\{80,{80},\overline{80},{70}\},\{{70},\overline{70},\overline{67},60\},
  \{{60},\overline{60},\overline{55},\overline{53}\},\{\overline{55},\overline{53},\overline{50},\overline{47}\},
   \{\overline{53},\overline{50},\overline{47},\overline{45}\},\]
   \[\{\overline{50},\overline{47},
   \overline{45},\overline{43}\},\{\overline{27},{20},{20},\overline{20}\},\{{20},{20},\overline{20},\overline{13}\},
   \{\overline{13},\overline{10},\overline{7},\overline{5}\},\{{ \overline{10}},\overline{7},{ \overline{5}},\overline{3}\}.\]
It can be checked that the first six $4$-bands are even and the last four $4$-bands are odd.

The following lemma says that two $(k-1)$-bands with overlapping parts are of the same parity.

\begin{lem}\label{parity-k-1-sequence-over-old} Let $\pi=(\pi_1,\pi_2,\ldots, \pi_\ell)$ be an overpartition   in ${\mathcal{\overline{B}}}(\alpha_1,\ldots,\alpha_\lambda;\eta,k,r)$ and let  $\{\pi_{c+l}\}_{0\leq l\leq k-2}$ and $\{\pi_{d+l}\}_{0\leq l\leq k-2}$   be two $(k-1)$-bands  of $\pi$ with $c< d\leq c+k-2$. Then, $\{\pi_{c+l}\}_{0\leq l\leq k-2}$ and $\{\pi_{d+l}\}_{0\leq l\leq k-2}$  are of the same parity.
  \end{lem}

\pf To show that $\{\pi_{c+l}\}_{0\leq l\leq k-2}$ and $\{\pi_{d+l}\}_{0\leq l\leq k-2}$ are of the same parity, it is equivalent to show that
\begin{equation}\label{parity-equiva-0}
\begin{split}
&\quad\left[|\pi_c|/\eta\right]+\cdots+\left[|\pi_{c+k-2}|/\eta\right]+\overline{V}_\pi(\pi_c)+\overline{O}_\pi(\pi_{c+k-2})\\[5pt]
&\equiv\left[|\pi_d|/\eta\right]+\cdots+\left[|\pi_{d+k-2}|/\eta\right]+\overline{V}_\pi(\pi_d)+\overline{O}_\pi(\pi_{d+k-2})\pmod{2}.
\end{split}
\end{equation}

Assume that $d=c+t$ where $1\leq t\leq k-2$. The overlapping structure of  $\{\pi_{c+l}\}_{0\leq l\leq k-2}$ and  $\{\pi_{d+l}\}_{0\leq l\leq k-2}$ can be described as follows:
{\footnotesize\[\begin{array}{ccccccccccccccccccc}
\pi_{d+k-2}&\leq &\cdots
&\leq& \pi_{d+k-1-t}&<&\pi_{d+k-2-t}&\leq&\cdots&\leq& \pi_d\\[5pt]
&&&&&&\|&&&&\|&&&&\\[5pt]
&&&&&&\pi_{c+k-2} &\leq& \cdots &\leq &\pi_{c+t}&<&\pi_{c+t-1}&\leq
&\cdots&\leq&\pi_{c}.
\end{array}
\]}
Then, we have
\[\begin{split}&\quad\overline{V}_\pi(\pi_c)-\overline{V}_\pi(\pi_d)\\
&=\text{the number of overlined parts not divisible by }\eta\text{ in }\pi_{c},\ldots,\pi_{c+t-1},
\end{split}
\]
\[\begin{split}&\quad\overline{O}_\pi(\pi_{d+k-2})-\overline{O}_\pi(\pi_{c+k-2})\\
&=\text{the number of overlined parts divisible by }\eta\text{ in }\pi_{d+k-1-t},\ldots,\pi_{d+k-2},
\end{split}\]
and
 \begin{eqnarray*}
   &&\quad\left[|\pi_c|/\eta\right]+\cdots+\left[|\pi_{c+k-2}|/\eta\right]-\left(\left[|\pi_d|/\eta\right]+\cdots+\left[|\pi_{d+k-2}|/\eta\right]\right)\\[5pt]
  &&=
  [|\pi_{c}|/\eta]+\cdots
   +[|\pi_{c+t-1}|/\eta]-
   \left([|\pi_{d+k-1-t}|/\eta]+\cdots+[|\pi_{d+k-2}|/\eta]\right).
   \end{eqnarray*}

   We find that in order to show \eqref{parity-equiva-0}, it suffices to prove that
   \begin{equation} \label{cong-g-x-f-over}
   \begin{split}
   &\quad[|\pi_{c+t-1}|/\eta]+\cdots
   +[|\pi_{c}|/\eta]-
   \left([|\pi_{d+k-2}|/\eta]+\cdots+[|\pi_{d+k-1-t}|/\eta]\right)\\[5pt]
   &\equiv\text{the number of overlined parts divisible by }\eta\text{ in }\pi_{d+k-1-t},\ldots,\pi_{d+k-2} \\[5pt]
   &\quad+ \text{the number of overlined parts not divisible by }\eta\text{ in }\pi_{c},\ldots,\pi_{c+t-1}\pmod{2}.
  \end{split}
   \end{equation}

   We consider the following two cases.

   Case 1: $\pi_c$ is divisible by $\eta$. In this case, we may write $\pi_c=(b+2)\eta$ or $\overline{(b+2)\eta}$. Then, we have
   \[b\eta \leq \pi_{d+k-2}\leq   \cdots \leq \pi_{d+k-1-t}<(b+1)\eta,\]
   and
\[(b+1)\eta<\pi_{c+t-1}\leq  \cdots \leq \pi_{c}\leq(b+2) \eta.\]
This implies that for $k-1-t\leq l\leq k-2$,
$[|\pi_{d+l}|/\eta]=b+1$ if $\pi_{d+l}$ is an overlined part divisible by $\eta$, or otherwise $[|\pi_{d+l}|/\eta]=b$,  and  for  $0\leq l\leq t-1$, $[\pi_{c+l}/\eta]=b+1$  if $\pi_{c+l}$ is an overlined part not divisible by $\eta$, or otherwise $[\pi_{c+l}/\eta]=b+2$.  So, \eqref{cong-g-x-f-over} is confirmed.

Case 2: $\pi_c$ is not divisible by $\eta$. In this case, we may write  $\pi_c=\overline{(b+1)\eta+\alpha_s}$, where $1\leq s\leq \lambda$. Then, we have
\begin{equation}\label{lemma2.6-pfa}
     \overline{(b-1)\eta+\alpha_s}<\pi_{d+k-2}\leq   \cdots \leq \pi_{d+k-1-t}\leq \overline{b\eta+\alpha_s},
 \end{equation}
  and
  \begin{equation}\label{lemma2.6-pfb}
      \overline{b\eta+\alpha_s}<\pi_{c+t-1}\leq   \cdots \leq \pi_{c}=\overline{(b+1)\eta+\alpha_s}.
  \end{equation}
  Assume that there are $f_1$ parts $\pi_{d+l}$ in \eqref{lemma2.6-pfa} satisfying $\overline{(b-1)\eta+\alpha_{s}}< \pi_{d+l}\leq\overline{(b-1)\eta+\alpha_\lambda}.$
 For such  a part $\pi_{d+l}$, we have $[|\pi_{d+l}|/\eta]=b-1$. Assume that there are $f_2$ parts $\pi_{d+l}$ in \eqref{lemma2.6-pfa} satisfying $\pi_{d+l}=\overline{b\eta}.$ In this case, we have $[|\pi_{d+l}|/\eta]=b$. Assume that there are $f_3$ parts $\pi_{d+l}$  in \eqref{lemma2.6-pfa} satisfying
$b\eta\leq \pi_{d+l}\leq \overline{b\eta+\alpha_s},$  which gives $[|\pi_{d+l}|/\eta]=b$. Then, there are $f_2$ overlined parts divisible by $\eta$ in $\pi_{d+k-1-t},\ldots,\pi_{d+k-2}$. Clearly, we have
\begin{equation}\label{f-equal-1}
f_1+f_2+f_3=t.
\end{equation}

 Assume that there are $f_4$ parts $\pi_{c+l}$  in \eqref{lemma2.6-pfb} satisfying
$\overline{b\eta+\alpha_{s}}< \pi_{c+l}\leq \overline{b\eta+\alpha_\lambda}.$ For such  a part $\pi_{c+l}$, we have $[|\pi_{c+l}|/\eta]=b$. Assume that there are $f_5$ parts $\pi_{c+l}$   in \eqref{lemma2.6-pfb} satisfying
  $\pi_{c+l}=\overline{(b+1)\eta}$ or $(b+1)\eta$. In this case, we have   $[|\pi_{c+l}|/\eta]=b+1$. Assume that there are $f_6$ parts
  $\pi_{c+l}$ in \eqref{lemma2.6-pfb} satisfying
  $\overline{(b+1)\eta+\alpha_1}\leq \pi_{c+l}\leq \overline{(b+1)\eta+\alpha_s},$ which gives  $[\pi_{c+l}/\eta]=b+1$. Then, there are $f_4+f_6$ overlined parts not divisible by $\eta$ in $\pi_{c},\ldots,\pi_{c+t-1}$. Clearly, we have
\begin{equation}\label{f-equal-2}
f_4+f_5+f_6=t.
\end{equation}

We proceed to show that
\begin{equation}\label{f-equal-3}
f_1+f_2=f_4+f_5\text{ and }f_3=f_6.
\end{equation}
By virtue of \eqref{f-equal-1} and \eqref{f-equal-2}, it suffices to prove that $f_3=f_6$. Since $\pi$ satisfies the condition (3) in Definition \ref{defi-O-B}, we have $f_3+k-t-1+f_4+f_5\leq k-1$, that is, $f_3+f_4+f_5\leq t$. Combining with \eqref{f-equal-2},  we get $f_3\leq f_6$. So, in order to show that $f_3=f_6$, it remains to verify that $f_3\geq f_6$. There are two cases.

(1) $\pi_d\geq\overline{(b+1)\eta}$. In such case, we have $\pi_{d+k-2}\geq b\eta$, and so
\[{b\eta}\leq\pi_{d+k-2}\leq   \cdots \leq \pi_{d+k-1-t}\leq \overline{b\eta+\alpha_s}.
\]
It yields that $f_3=t\geq f_6$.

(2) $ \overline{b\eta+\alpha_s}<\pi_d<\overline{(b+1)\eta}$. In this case, we may write $\pi_d=\overline{b\eta+\alpha_g}$ with $s<g\leq \lambda$. In this case, we have $\pi_{d+k-2}>\overline{(b-1)\eta+\alpha_g}$, and so $f_1$ is the number of parts $\pi_{d+l}$ in \eqref{lemma2.6-pfa} satisfying $\overline{(b-1)\eta+\alpha_{g}}< \pi_{d+l}\leq\overline{(b-1)\eta+\alpha_\lambda}.$

Under the assumption that $\pi_c=\overline{(b+1)\eta+\alpha_s}$, we get $\pi_{d+k-2-t}=\pi_{c+k-2}>\pi_c-\eta=\overline{b\eta+\alpha_s}$. So, we have
\[\overline{b\eta+\alpha_s}<\pi_{d+k-2-t}\leq   \cdots \leq \pi_{d}=\overline{b\eta+\alpha_g}.
\]
It follows that $k-t-1\leq g-s$, that is, $k-1-g+s\leq t$. Given the condition that $k-1>\lambda$, we obtain that
\[f_1+f_2+f_6\leq (\lambda-g)+1+s\leq k-1-g+s\leq t.\]
Combining with \eqref{f-equal-1}, we get $f_3\geq f_6$.

Now, we conclude that $f_3\geq f_6$, and so \eqref{f-equal-3} is satisfied.  It yields that
\begin{equation*}\label{f-equal-4}
f_1+f_5\equiv f_2+f_4\pmod{2}.
\end{equation*}
Combining with \eqref{f-equal-1} and \eqref{f-equal-2}, we get
\[
\begin{split}
&\quad[|\pi_{c+t-1}|/\eta]+\cdots+[|\pi_{c}|/\eta]-\left([|\pi_{d+k-2}|/\eta]+\cdots+[|\pi_{d+k-1-t}|/\eta]\right)\\[5pt]
&=(b+1)(f_6+f_5)+bf_4-b(f_3+f_2)-(b-1)f_1\\[5pt]
&=f_6+f_5+f_1\\[5pt]
&\equiv f_6+f_4+f_2\pmod2.
\end{split}
\]

Recall that the number of overlined parts divisible by $\eta$  in $\pi_{d+k-1-t},\ldots,\pi_{d+k-2}$  is $f_2$ and the number of overlined parts not divisible by $\eta$ in $\pi_{c},\ldots,\pi_{c+t-1}$ is $f_4+f_6$, so \eqref{cong-g-x-f-over} is valid. This completes the proof.    \qed

We conclude this section with the following lemma.

\begin{lem}\label{cong-t-0}
Let $\pi$ be an overpartition in $\overline{\mathcal{B}}(\alpha_1,\ldots,\alpha_\lambda;\eta,k,r)$. For $t\geq 1$, assume that $\{\pi_{c+l}\}_{0\leq l\leq k-2}$ and $\{\pi_{d+l}\}_{0\leq l\leq k-2}$  are two $(k-1)$-bands  of $\pi$ belonging to $[(t-1)\eta,(t+1)\eta]$.
Then,  $\{\pi_{c+l}\}_{0\leq l\leq k-2}$ and $\{\pi_{d+l}\}_{0\leq l\leq k-2}$ are of the same parity.
\end{lem}

\pf We consider the following two cases.

Case 1: $\{\pi_{c+l}\}_{0\leq l\leq k-2}$ and $\{\pi_{d+l}\}_{0\leq l\leq k-2}$ have overlapping parts. It is a consequence of  Lemma \ref{parity-k-1-sequence-over-old}.

Case 2: $\{\pi_{c+l}\}_{0\leq l\leq k-2}$ and $\{\pi_{d+l}\}_{0\leq l\leq k-2}$ do not have overlapping parts. Without loss of generality, we assume that $\pi_c>\pi_d$. Then, we have
\[(t-1)\eta\leq \pi_{d+k-2}\leq\cdots\leq \pi_d<\pi_{c+k-2}\leq\cdots\leq \pi_c\leq (t+1)\eta.\]

Under the condition that $\pi\in\overline{\mathcal{B}}(\alpha_1,\ldots,\alpha_\lambda;\eta,k,r)$, we see that there are at most $k-1$ of the $\pi_i$ such that $(t-1)\eta\leq \pi_i\leq t\eta$ and there are at most $k-1$ of the $\pi_i$ such that $t\eta\leq \pi_i\leq (t+1)\eta$. So, there are at most $2k-2$ of the $\pi_i$ such that $(t-1)\eta\leq \pi_i\leq (t+1)\eta$. Therefore, $\pi_{d+k-2},\ldots,\pi_d,\pi_{c+k-2},\ldots,\pi_c$ are the $2k-2$ of the $\pi_i$ such that $(t-1)\eta\leq \pi_i\leq (t+1)\eta$. Moreover, we have
\begin{equation}\label{d-parts}
(t-1)\eta\leq \pi_{d+k-2}\leq\cdots\leq \pi_d<t\eta,
\end{equation}
and
\begin{equation}\label{c-parts}
t\eta<\pi_{c+k-2}\leq\cdots\leq \pi_c\leq (t+1)\eta.
\end{equation}

 Assume that there are $f_1$ parts $\pi_{d+l}$ in \eqref{d-parts} satisfying ${(t-1)\eta}\leq \pi_{d+l}\leq\overline{(t-1)\eta+\alpha_\lambda}.$
 For such  a part $\pi_{d+l}$, we have $[|\pi_{d+l}|/\eta]=t-1$. Assume that there are $f_2$ parts $\pi_{d+l}$ in \eqref{d-parts} satisfying $\pi_{d+l}=\overline{t\eta},$ which gives $[|\pi_{d+l}|/\eta]=t$. Clearly,
 \begin{equation}\label{d-c-con-1}
 f_1+f_2=k-1\text{ and }\overline{O}_\pi(\pi_{d+k-2})-\overline{O}_\pi(\pi_{c+k-2})=f_2.
 \end{equation}

Assume that there are $f_3$ parts $\pi_{c+l}$ in \eqref{c-parts} satisfying $\overline{t\eta+\alpha_1}\leq \pi_{c+l}\leq\overline{t\eta+\alpha_\lambda}.$ In this case, we have $[|\pi_{c+l}|/\eta]=t$. Assume that there are $f_4$ parts $\pi_{c+l}$ in \eqref{c-parts} satisfying $\pi_{c+l}=\overline{(t+1)\eta}$ or $(t+1)\eta$, which implies that $[|\pi_{c+l}|/\eta]=t+1$.
 Clearly,
 \begin{equation*}\label{d-c-con-2}
 f_3+f_4=k-1\text{ and }\overline{V}_\pi(\pi_c)-\overline{V}_\pi(\pi_d)=f_3.
 \end{equation*}
Combining with \eqref{d-c-con-1}, we get
\begin{equation*}\label{equiv-c-d}
\begin{split}
&\quad\left[|\pi_c|/\eta\right]+\cdots+\left[|\pi_{c+k-2}|/\eta\right]-\left[|\pi_d|/\eta\right]+\cdots+\left[|\pi_{d+k-2}|/\eta\right]\\[5pt]
&=tf_3+(t+1)f_4-(t-1)f_1-tf_2\\[5pt]
&=f_3+2f_4-f_2\\[5pt]
&\equiv f_3+f_2\\[5pt]
&= \overline{V}_\pi(\pi_c)-\overline{V}_\pi(\pi_d)+\overline{O}_\pi(\pi_{d+k-2})-\overline{O}_\pi(\pi_{c+k-2})\pmod2,
\end{split}
\end{equation*}
which implies that $\{\pi_{c+l}\}_{0\leq l\leq k-2}$ and $\{\pi_{d+l}\}_{0\leq l\leq k-2}$ are of the same parity. The proof is completed.  \qed

\section{Proof of Theorem \ref{eqv-main-0} }

The goal of this section  is to give a proof of Theorem \ref{eqv-main-0}. To do this,  we shall define the  $(k-1)$-reduction operation  and   the $(k-1)$-augmentation operation, which are the main ingredients in the construction of $\Phi$ in Theorem \ref{eqv-main-0}.

\subsection{The $(k-1)$-reduction and the $(k-1)$-augmentation}

 The definitions of the  $(k-1)$-reduction and the $(k-1)$-augmentation are based on two subsets of $\overline{\mathcal{B}}_0(\alpha_1,\ldots,\alpha_\lambda;\eta,k,r)$. To describe these two subsets, we need
to introduce the following notation.  Define $s(\pi)$ to be the smallest overlined part divisible by $\eta$ in $\pi$ with the convention that $s(\pi)=+\infty$ if there are no   overlined parts  divisible by $\eta$ in $\pi$. Define $g(\pi)$ to be the smallest part $\pi_i$ such that  $\{\pi_{i+l}\}_{0\leq l\leq k-2}$ is a $(k-1)$-band of $\pi$  with the convention that  $g(\pi)=+\infty$ if there  are no $(k-1)$-bands of   $\pi$.

Let $\pi$ be an overpartition in $\overline{\mathcal{B}}_0(\alpha_1,\ldots,\alpha_\lambda;\eta,k,r)$. Under the condition that $k-1>\lambda$, we derive that $g(\pi)\geq \overline{\eta}$. Moreover, we have $g(\pi)\geq \eta$ if $s(\pi)>\overline{\eta}$.
So, the condition (5) in Definition \ref{defi-O-B} can be expressed as follows:
\[\text{if }f_{\leq \eta}(\pi)=r\text{ and }s(\pi)>\overline{\eta},\text{ then }\eta\leq g(\pi)<\overline{2\eta}.\]
Consequently, we have $f_{\leq \eta}(\pi)\leq r-1$ if $s(\pi)>\overline{\eta}$ and $g(\pi)\geq\overline{2\eta}$.

 We are now in a position to present two subsets of $\overline{\mathcal{B}}_0(\alpha_1,\ldots,\alpha_\lambda;\eta,k,r)$.

\begin{itemize}

\item  For $t\geq 1$, let $\overline{\mathcal{B}}^{\,=}_0(\alpha_1,\ldots,\alpha_\lambda;\eta,k,r|t)$ denote the set of   overpartitions $\pi$ in $\overline{\mathcal{B}}_0(\alpha_1,\ldots,\break $$\alpha_\lambda;\eta,k,r)$ such that either $s(\pi)=\overline{t\eta}$ and $g(\pi)\geq\overline{t\eta}$, or $s(\pi)>\overline{t\eta}$ and  ${t\eta}\leq g(\pi)<\overline{(t+1)\eta}$.

\item For $t\geq 1$, let $\overline{\mathcal{B}}^{\,>}_0(\alpha_1,\ldots,\alpha_\lambda;\eta,k,r|t)$ denote the set of   overpartitions $\pi$ in $\overline{\mathcal{B}}_0(\alpha_1,\ldots,\break$$\alpha_\lambda;\eta,k,r)$ such that
 $s(\pi)>\overline{t\eta}$ and $g(\pi)\geq \overline{(t+1)\eta}$.
\end{itemize}

Let $\pi$ be an overpartition in $\overline{\mathcal{B}}^{\,=}_0(\alpha_1,\ldots,\alpha_\lambda;\eta,k,r|t)$ with $t\geq 2$ or $\overline{\mathcal{B}}^{\,>}_0(\alpha_1,\ldots,\alpha_\lambda;\eta,k,r|t)$ with $t\geq 1$. Then, we have $s(\pi)>\overline{\eta}$ and $g(\pi)\geq\overline{2\eta}$, and so $f_{\leq \eta}(\pi)\leq r-1$.

The following proposition provides a criterion to determine whether an overpartition in $\overline{\mathcal{B}}^{\,=}_0(\alpha_1,\ldots,\alpha_\lambda;\eta,k,r|t)$ is also an overpartition in $\overline{\mathcal{B}}^{\,>} _0(\alpha_1,\ldots,\alpha_\lambda;\eta,k,r|t')$.

\begin{prop}\label{daxiao-0-4}
For $t\geq 2$, let $\pi$ be an overpartition in $\overline{\mathcal{B}}^{\,=}_0(\alpha_1,\ldots,\alpha_\lambda;\eta,k,r|t)$. Then $\pi$ is an overpartition in $\overline{\mathcal{B}}^{\,>} _0(\alpha_1,\ldots,\alpha_\lambda;\eta,k,r|t')$ if and only if $t'<t$.
\end{prop}

\pf  By definition, we see that $\pi$ is  an overpartition in $\overline{\mathcal{B}}^{\,=}_0(\alpha_1,\ldots,\alpha_\lambda;\eta,k,r|t)$ if and only if
$\pi$ is  an overpartition in $\overline{\mathcal{B}}_0(\alpha_1,\ldots,\alpha_\lambda;\eta,k,r)$ such that
\begin{equation}\label{daxiao-eqn-0-4}
 \min\left\{\left[ |s(\pi)|/\eta\right],\left[|g(\pi)|/\eta\right]\right\}=t.
    \end{equation}

On the other hand,   $\pi$ is  an overpartition in $\overline{\mathcal{B}}^{\, >}_0(\alpha_1,\ldots,\alpha_\lambda;\eta,k,r|t')$ if and only if $\pi$ is  an overpartition in $\overline{\mathcal{B}}_0(\alpha_1,\ldots,\alpha_\lambda;\eta,k,r)$ such that
\begin{equation}\label{daxiao-eqn-0-4b}
 \min\left\{\left[ |s(\pi)|/\eta\right],\left[ |g(\pi)|/\eta\right]\right\}\geq t'+1.
    \end{equation}
 Combining  \eqref{daxiao-eqn-0-4} and \eqref{daxiao-eqn-0-4b} completes the proof.  \qed

The $(k-1)$-reduction operation is stated as follows.

\begin{defi}[The $(k-1)$-reduction]\label{defi-division} For  $t\geq 1$,   let $\pi$ be an overpartition in $\overline{\mathcal{B}}^{\,=}_0(\alpha_1,\ldots,\break$
$\alpha_\lambda;\eta,k,r|t)$. Define the $(k-1)$-reduction $D_t\colon \pi \rightarrow \mu$ as follows{\rm{:}}    If $s(\pi)=\overline{t\eta}$, then $\mu$ is obtained  from $\pi$ by removing the overlined part $\overline{t\eta}$. Otherwise, $\mu$ is obtained  from $\pi$ by removing a non-overlined part   ${t\eta}$.
\end{defi}

The following proposition ensures that the $(k-1)$-reduction is well defined.

\begin{prop}\label{daxiao-0-5}
For $t\geq1$, let $\pi$ be an overpartition in $\overline{\mathcal{B}}_0(\alpha_1,\ldots,\alpha_\lambda;\eta,k,r)$ such that $s(\pi)>\overline{t\eta}$ and  ${t\eta}\leq g(\pi)<\overline{(t+1)\eta}$. Then,  $t\eta$ occurs in $\pi$.
\end{prop}

\pf Suppose to the contrary that $t\eta$ does not occur in $\pi$. In this case, we have $g(\pi)\neq t\eta$. Under that condition that ${t\eta}\leq g(\pi)<\overline{(t+1)\eta}$, we have ${t\eta}< g(\pi)<\overline{(t+1)\eta}$. So, we can write $g(\pi)=\overline{t\eta+\alpha_s}$, where $1\leq s\leq \lambda$.

Assume that $g(\pi)$ is the $i$-th part of $\pi$,  then we have $\pi_i=g(\pi)=\overline{t\eta+\alpha_s}$ and $\{\pi_{i+l}\}_{0\leq l\leq k-2}$ is a $(k-1)$-band of $\pi$. So, we get
\[\overline{(t-1)\eta+\alpha_s}<\pi_{i+k-2}\leq \cdots\leq \pi_i=\overline{t\eta+\alpha_s}.\]

For $0\leq l\leq k-2$, it follows from $s(\pi)>\overline{t\eta}$ that $\overline{t\eta}$ does not occur in $\pi$, and so  $\pi_{i+l}\neq \overline{t\eta}$. Under the assumption that  $t\eta$ does not occur in $\pi$, we have $\pi_{i+l}\neq t\eta$. So, we obtain that
\[\overline{(t-1)\eta+\alpha_s}<\pi_{i+l}\leq \overline{(t-1)\eta+\alpha_\lambda}\quad\text{or}\quad\overline{t\eta+\alpha_1}\leq\pi_{i+l}\leq \overline{t\eta+\alpha_s}.\]
It implies that $k-1\leq (\lambda-s)+s=\lambda$, which contradicts the fact that $k-1>\lambda$. This completes the proof.   \qed

The following lemma tell us that the $(k-1)$-reduction operation is a map from $\overline{\mathcal{B}}^{\,=}_0(\alpha_1,\ldots,\alpha_\lambda;\eta,k, r|t)$ to $\overline{\mathcal{B}}^{\,>}_0(\alpha_1,\ldots,\alpha_\lambda;\eta,k, r|t)$.

 \begin{lem}\label{division} For  $t\geq 1$,  let $\pi$ be an overpartition in $\overline{\mathcal{B}}^{\,=}_0(\alpha_1,\ldots,\alpha_\lambda;\eta,k, r|t)$ and let $\mu=D_t(\pi)$. Then $\mu$ is overpartition in $\overline{\mathcal{B}}^{\,>}_0(\alpha_1,\ldots,\alpha_\lambda;\eta,k,r|t)$. Furthermore, $|\mu|=|\pi|-t\eta$ and $\ell(\mu)=\ell(\pi)-1$.
\end{lem}

\pf By definition, we wish to show that $\mu$ satisfies the following conditions:
\begin{itemize}
    \item[(A)] $\mu$ is an overpartition in $\overline{\mathcal{B}}(\alpha_1,\ldots,\alpha_\lambda;\eta,k,r)$;

    \item [(B)] $s(\mu)>\overline{t\eta}$ and $g(\mu)\geq \overline{(t+1)\eta}$;

    \item [(C)] All $(k-1)$-bands  of $\mu$  are even.
\end{itemize}

{\noindent Condition (A).} Given  the precondition  $\pi \in \overline{\mathcal{B}}^{\,=}_0(\alpha_1,\ldots,\alpha_\lambda;\eta,k, r|t)$, it is immediate from the construction of $\mu$   that it satisfies the conditions (1)-(4) in Definition \ref{defi-O-B}. 
That is to say,  $\mu$ is an overpartition in $\overline{\mathcal{B}}(\alpha_1,\ldots,\alpha_\lambda;\eta,k,r)$.

{\noindent Condition (B).} Since $\pi \in \overline{\mathcal{B}}^{\,=}_0(\alpha_1,\ldots,\alpha_\lambda;\eta,k, r|t)$, we have $s(\pi)\geq\overline{t\eta}$ and $g(\pi)\geq \overline{t\eta}$. From the construction of $\mu$, we deduce that $s(\mu)\geq s(\pi)\geq\overline{t\eta}$ and $g(\mu)\geq g(\pi)\geq \overline{t\eta}$. Moreover, we see that $\overline{t\eta}$ does not occur in $\mu$, which implies that  $s(\mu)>\overline{t\eta}$. It follows from Lemmas \ref{t+1} and \ref{t+1-o} that there are no $(k-1)$-bands of $\mu$ belonging to $[(t-1)\eta,\overline{(t+1)\eta})$, and so $g(\mu)\geq \overline{(t+1)\eta}$. Thus, the condition (B) holds.

{\noindent Condition (C).} Assume that
 $\{\mu_{i+l}\}_{0\leq l\leq k-2}$ is a $(k-1)$-band of $\mu$. We aim to show that $\{\mu_{i+l}\}_{0\leq l\leq k-2}$ is  even.  Using the condition (B), we know that $g(\mu)\geq \overline{(t+1)\eta}$, and so $\mu_{i}\geq \overline{(t+1)\eta}$. The assumption that $\{\mu_{i+l}\}_{0\leq l\leq k-2}$ is a $(k-1)$-band of $\mu$ yields $\mu_{i+k-2}\geq {t\eta}$. By the construction of $\mu$, we obtain that $\pi_{i+l}=\mu_{i+l}$ for $0\leq l\leq k-2$. Moreover, we have
 \begin{equation}\label{even-condition-0}
 \overline{V}_{\pi}(\pi_{i})=\overline{V}_{\mu}(\mu_{i})\text{ and }\overline{O}_\pi(\pi_{i+k-2})=\overline{O}_\mu(\mu_{i+k-2}).
 \end{equation}

  Clearly,  $\{\pi_{i+l}\}_{0\leq l\leq k-2}$ is a $(k-1)$-band of $\pi$. Since $\pi \in \overline{\mathcal{B}}^{\,=}_0(\alpha_1,\ldots,\alpha_\lambda;\eta,k, r|t)$, we have
   \begin{equation*}\label{div-pro}
 \left[|\pi_{i}|/\eta\right]+\left[|\pi_{i+1}|/\eta\right]+\cdots+\left[|\pi_{i+k-2}|/\eta\right]\equiv r-1+\overline{V}_{\pi}(\pi_{i})+\overline{O}_{\pi}(\pi_{i+k-2})\pmod{2}.
 \end{equation*}
  Combining with  \eqref{even-condition-0}, we obtain that $\{\mu_{i+l}\}_{0\leq l\leq k-2}$ of $\mu$ is even, and this proves that the condition (C) is justified.

In conclusion,  we have shown that   $\mu$ is an overpartition in  $\overline{\mathcal{B}}^{\,>}_0(\alpha_1,\ldots,\alpha_\lambda;\eta,k,r|t)$.  Clearly,   $|\mu|=|\pi|-t\eta$ and $\ell(\mu)=\ell(\pi)-1$. This completes the proof.  \qed

Before giving the definition of the $(k-1)$-augmentation operation, for $t\geq1$, we introduce the types of   $(k-2)$-bands of $\mu$ belonging to $[(t-1)\eta,\overline{(t+1)\eta})$, where $\mu$ is an overpartition in $\overline{\mathcal{B}}^{\,>}_0(\alpha_1,\ldots,\alpha_\lambda;\eta,k,r|t)$.  Assume that    $\{\mu_{i+l}\}_{0\leq l\leq k-3}$ is a  $(k-2)$-band of $\mu$ belonging to $[(t-1)\eta,\overline{(t+1)\eta})$, if
\[[|\mu_{i}|/\eta]+\cdots+[|\mu_{i+k-3}|/\eta]\equiv t+r-1+\overline{V}_{\mu}(\mu_{i})+\overline{O}_{\mu}(\mu_{i+k-3})\pmod2,\]
then we say that $\{\mu_{i+l}\}_{0\leq l\leq k-3}$ is of type $N$. Otherwise, we say that it is of type $O$.

\begin{prop}\label{k-1-band-con}
For $t\geq1$, let $\mu$ be an overpartition in $\overline{\mathcal{B}}^{\,>}_0(\alpha_1,\ldots,\alpha_\lambda;\eta,k,r|t)$ such that there exist $(k-2)$-bands of $\mu$ belonging to $[(t-1)\eta,\overline{(t+1)\eta})$. Then,
 \begin{itemize}
\item[{\rm(1)}]  all $(k-2)$-bands of $\mu$ belonging to $[(t-1)\eta,\overline{(t+1)\eta})$ are of the same type{\rm;}

\item[{\rm (2)}]  all $(k-2)$-bands of $\mu$ belonging to $[(t-1)\eta,\overline{(t+1)\eta})$ are of type O if $g(\mu)=\overline{(t+1)\eta}$ or $(t+1)\eta$.

\end{itemize}
\end{prop}

\pf Let $\pi$ be the overpartition obtained by inserting $\overline{t\eta}$ as an overlined part into $\mu$. Under the condition that $\mu\in\overline{\mathcal{B}}^{\,>}_0(\alpha_1,\ldots,\alpha_\lambda;\eta,k,r|t)$, we have $g(\mu)\geq \overline{(t+1)\eta}$ and  $f_{\leq \eta}(\mu)\leq r-1$. Clearly, $\pi$ is an overpartition in $\overline{\mathcal{B}}(\alpha_1,\ldots,\alpha_\lambda;\eta,k,r)$.

 Assume that $\{\mu_{m+l}\}_{0\leq l\leq k-3}$ is a $(k-2)$-band of $\mu$ belonging to $[(t-1)\eta,\overline{(t+1)\eta})$. It follows from the condition (2) in Lemma \ref{t+1-o} that $\{\pi_{m+l}\}_{0\leq l\leq k-2}$ is a $(k-1)$-band of $\pi$ belonging to $[(t-1)\eta,\overline{(t+1)\eta})$.  Assume that
\[[|\pi_{m}|/\eta]+\cdots+[|\pi_{m+k-2}|/\eta]\equiv a+\overline{V}_{\pi}(\pi_{m})+\overline{O}_{\pi}(\pi_{m+k-2})\pmod2.\]

Appealing to \eqref{cong-m-band-o} in Lemma \ref{t+1-o}, we get
\begin{equation*}\label{k-2-band-congruence}
[|\mu_{m}|/\eta]+\cdots+[|\mu_{m+k-3}|/\eta]
\equiv a+t+1+\overline{V}_{\mu}(\mu_{m})+\overline{O}_{\mu}(\mu_{m+k-3})\pmod2.
\end{equation*}

{\noindent Condition (1).} We are obliged  to show that for any $(k-2)$-band $\{\mu_{i+l}\}_{0\leq l\leq k-3}$  of $\mu$ belonging to $[(t-1)\eta,\overline{(t+1)\eta})$, $\{\mu_{i+l}\}_{0\leq l\leq k-3}$ and $\{\mu_{m+l}\}_{0\leq l\leq k-3}$ are of the same type, that is,
\begin{equation}\label{k-2-band-congruence-(1)}
[|\mu_{i}|/\eta]+\cdots+[|\mu_{i+k-3}|/\eta]
\equiv a+t+1+\overline{V}_{\mu}(\mu_{i})+\overline{O}_{\mu}(\mu_{i+k-3})\pmod2.
\end{equation}

  Again by Lemma \ref{t+1-o}, we see that $\{\pi_{i+l}\}_{0\leq l\leq k-2}$ is a $(k-1)$-band of $\pi$. Moreover, $\{\pi_{i+l}\}_{0\leq l\leq k-2}$ and  $\{\pi_{m+l}\}_{0\leq l\leq k-2}$ have overlapping part $\overline{t\eta}$. Utilizing Lemma \ref{parity-k-1-sequence-over-old}, we derive that $\{\pi_{i+l}\}_{0\leq l\leq k-2}$ and  $\{\pi_{m+l}\}_{0\leq l\leq k-2}$ are of same parity, and so
\[[|\pi_{i}|/\eta]+\cdots+[|\pi_{i+k-2}|/\eta]\equiv a+\overline{V}_{\pi}(\pi_{i})+\overline{O}_{\pi}(\pi_{i+k-2})\pmod2.\]
Using \eqref{cong-m-band-o} in Lemma \ref{t+1-o}, we can prove that \eqref{k-2-band-congruence-(1)} holds.  So, the condition (1) is verified.

{\noindent Condition (2).} Assume that $g(\mu)$ is the $s$-th part of $\mu$.  Then, we have $\mu_s=\overline{(t+1)\eta}$ or $(t+1)\eta$. Moreover, $\{\mu_{s+l}\}_{0\leq l\leq k-2}$ is a $(k-1)$-band of $\mu$, and so $\mu_{s+k-2}\geq t\eta$. Since $\mu\in\overline{\mathcal{B}}^{\,>}_0(\alpha_1,\ldots,\alpha_\lambda;\eta,k,r|t)$, we see that the $(k-1)$-band $\{\mu_{s+l}\}_{0\leq l\leq k-2}$ of $\mu$ is even.

It follows from the construction of $\pi$ that $\pi_{s+l}=\mu_{s+l}$ for $0\leq l\leq k-2$. Then, we have
$t\eta\leq \pi_{s+k-2}\leq\cdots\leq \pi_{s}\leq (t+1)\eta.$
Moreover,  $\{\pi_{s+l}\}_{0\leq l\leq k-2}$ is a $(k-1)$-band of $\pi$ and it is even.

By virtue of Lemma \ref{cong-t-0}, we see that the $(k-1)$-band $\{\pi_{m+l}\}_{0\leq l\leq k-2}$ of $\pi$ is even, that is, $a\equiv r-1\pmod 2$. Combining with \eqref{k-2-band-congruence-(1)}, we can derive that the $(k-2)$-band $\{\mu_{m+l}\}_{0\leq l\leq k-3}$ of $\mu$ is of type O. Under the condition (1), we obtain that all $(k-2)$-bands of $\mu$ belonging to $[(t-1)\eta,\overline{(t+1)\eta})$ are of type O. This completes the proof.
 \qed

We are now in a position to introduce the $(k-1)$-augmentation operation.

\begin{defi}[The $(k-1)$-augmentation]\label{insertion} For $t\geq 1$, let $\mu$ be an overpartition in $\overline{\mathcal{B}}^{\,>}_0(\alpha_1,\ldots,\alpha_\lambda;\eta,k,r|t)$.
The $(k-1)$-augmentation  $C_t\colon \mu \rightarrow \pi$ is defined as follows{\rm{:}} If there exist  $(k-2)$-bands of $\mu$ belonging to $[(t-1)\eta,\overline{(t+1)\eta})$ which are of type N, then $\pi$ is obtained from $\mu$ by inserting a non-overlined part $t\eta$.  Otherwise,  $\pi$ is obtained from $\mu$ by inserting an overlined part $\overline{t\eta}$.

\end{defi}

The following lemma says that the $(k-1)$-augmentation operation is a map from $\overline{\mathcal{B}}^{\,>}_0(\alpha_1,\ldots,\alpha_\lambda;\eta,k, r|t)$ to $\overline{\mathcal{B}}^{\,=}_0(\alpha_1,\ldots,\alpha_\lambda;\eta,k, r|t)$.

\begin{lem}\label{combinatio}    For  $t\geq 1$,  let $\mu$ be an overpartition in $\overline{\mathcal{B}}^{\,>}_0(\alpha_1,\ldots,\alpha_\lambda;\eta,k, r|t)$ and let $\pi=C_t(\mu)$. Then   $\pi$ is an overpartition in $\overline{\mathcal{B}}^{\,=}_0(\alpha_1,\ldots,\alpha_\lambda;\eta,k,r|t)$ such that $|\pi|=|\mu|+t\eta$ and $\ell(\pi)=\ell(\mu)+1$.

\end{lem}

\pf  To prove that $\pi$ is an overpartition in $\overline{\mathcal{B}}^{\,=}_0(\alpha_1,\ldots,\alpha_\lambda;\eta,k,r|t)$, we need to verify that $\pi$ satisfies the following conditions:

\begin{itemize}

 \item[(A)] $\pi$ is an overpartition in $\overline{\mathcal{B}}(\alpha_1,\ldots,\alpha_\lambda;\eta,k,r)$;

\item[(B)] $s(\pi)=\overline{t\eta}$ and $g(\pi)\geq\overline{t\eta}$, or $s(\pi)>\overline{t\eta}$ and  ${t\eta}\leq g(\pi)<\overline{(t+1)\eta}$;

 \item[(C)]  If $f_{\leq \eta}(\pi)=r$ and $s(\pi)>\overline{\eta}$, then $\eta\leq g(\pi)<\overline{2\eta}$;

\item[(D)] All $(k-1)$-bands  of $\pi$  are even.
\end{itemize}

{\noindent Condition (A).} It is equivalent to show that $\pi$ satisfies the conditions (1)-(4) in Definition \ref{defi-O-B}. Clearly, $\pi$ satisfies the conditions (1) and (2) in Definition \ref{defi-O-B}. We proceed to show that $\pi$ satisfies the condition  (3)  in Definition \ref{defi-O-B}. Suppose to the contrary that   there exists  $1\leq c\leq \ell(\pi)-k+1$  such that
\begin{equation}\label{lema-ctd}
 \pi_c\leq \pi_{c+k-1}+\eta  \text{ with strict inequality if } \pi_c  \text{ is overlined.}
\end{equation}
 Assume that the $m$-th part of $\pi$ is the new inserted part of the $(k-1)$-augmentation operation, that is, $\pi_m=\overline{t\eta}$ or $t\eta$.  Since $\mu\in\overline{\mathcal{B}}^{\,>}_0(\alpha_1,\ldots,\alpha_\lambda;\eta,k, r|t)$, we have  $\pi_c\geq \pi_m\geq \pi_{c+k-1}$. Comparing with \eqref{lema-ctd}, we get $\pi_c\leq (t+1)\eta$ and $\pi_{c+k-1}\geq (t-1)\eta$.

  By the construction of $\pi$, we see that the consecutive $k$ parts $\pi_c,\ldots,\pi_{c+k-1}$ of $\pi$ consist of the consecutive $k-1$ parts $\mu_c,\ldots,\mu_{c+k-2}$ of $\mu$ together with the new inserted part $\overline{t\eta}$ or $t\eta$. Then, we arrive at
\[(t+1)\eta\geq\mu_c\geq \cdots \geq \mu_{c+k-2}\geq(t-1)\eta.\]
Theorefore, $\{\mu_{c+l}\}_{0\leq l\leq k-2}$ is a $(k-1)$-band of $\mu$ belonging to $[(t-1)\eta,{(t+1)\eta}]$.  The condition  $g(\mu)\geq\overline{(t+1)\eta}$ implies that $\mu_{c}=\overline{(t+1)\eta}$ or $(t+1)\eta$,  and so $\pi_{c}=\mu_{c}=\overline{(t+1)\eta}$ or $(t+1)\eta$. By \eqref{lema-ctd}, we get $\pi_{c+k-1}\geq t\eta$.
It yields that  $\pi_m=t\eta$, that is, $\pi$ is obtained from $\mu$ by inserting a non-overlined part $t\eta$. By definition, there exist  $(k-2)$-bands of $\mu$ belonging to $[(t-1)\eta,\overline{(t+1)\eta})$ which are of type N.
It follows from the condition (2) in Proposition \ref{k-1-band-con} that $g(\mu)\neq \overline{(t+1)\eta},(t+1)\eta$. It implies that $g(\mu)>(t+1)\eta$, which contradicts the condition that $\pi_{c}=\overline{(t+1)\eta}$ or $(t+1)\eta$. Thus, $\pi$ satisfies the condition  (3)  in Definition \ref{defi-O-B}.

By virtue of $\mu\in\overline{\mathcal{B}}^{\,>}_0(\alpha_1,\ldots,\alpha_\lambda;\eta,k, r|t)$, we have $f_{\leq \eta}(\mu)\leq r-1$. It follows from the construction of $\pi$ that  $f_{\leq \eta}(\pi)\leq f_{\leq \eta}(\mu)+1\leq r$. That is, $\pi$ satisfies the condition  (4)  in Definition \ref{defi-O-B}. Thus, the condition (A) is verified.

{\noindent Condition (B).}  Since $\mu\in\overline{\mathcal{B}}^{\,>}_0(\alpha_1,\ldots,\alpha_\lambda;\eta,k, r|t)$, we have $s(\mu)>\overline{t\eta}$ and $g(\mu)\geq \overline{(t+1)\eta}$. We consider the following two cases.

{\noindent Case (B)-1.}  $\pi$ is obtained from $\mu$ by inserting  an overlined part $\overline{t\eta}$. Obviously,   $s(\pi)=\overline{t\eta}$ and  $g(\pi)\geq \overline{t\eta}$.

{\noindent Case (B)-2.}  $\pi$ is obtained from $\mu$ by inserting a non-overlined part ${t\eta}$, that is, there exist  $(k-2)$-bands of $\mu$ belonging to $[(t-1)\eta,\overline{(t+1)\eta})$ which are of type N. Under the condition that $s(\mu)>\overline{t\eta}$ and $g(\mu)\geq \overline{(t+1)\eta}$, we deduce that $s(\pi)=s(\mu)>\overline{t\eta}$ and $g(\pi)\geq {t\eta}$.

Assume that $\{\mu_{d+l}\}_{0\leq l\leq k-3}$ is a $(k-2)$-band of $\mu$ belonging to $[(t-1)\eta,\overline{(t+1)\eta})$. It follows from Lemma \ref{t+1} that $\{\pi_{d+l}\}_{0\leq l\leq k-2}$ is a $(k-1)$-band of $\pi$ belonging to $[(t-1)\eta,\overline{(t+1)\eta})$. Moreover,
the consecutive $k-1$ parts $\pi_d,\ldots,\pi_{d+k-2}$ of $\pi$ consist of the parts $\mu_d,\ldots,\mu_{d+k-3}$ of $\mu$ together with the new inserted part ${t\eta}$. Hence, we have
\[(t-1)\eta\leq\pi_{d+k-2}\leq \cdots \leq \pi_{d}<\overline{(t+1)\eta},\]
which implies that  $\{\pi_{d+l}\}_{0\leq l\leq k-2}$ is a $(k-1)$-band of $\pi$ belonging to $[(t-1)\eta,\overline{(t+1)\eta})$. It yields that  $g(\pi)<\overline{(t+1)\eta}$. This
proves that the condition (B) is valid.

{\noindent Condition (C).} Assume that $f_{\leq \eta}(\pi)=r$ and $s(\pi)>\overline{\eta}$. In this case,  $\pi$ is obtained from $\mu$ by inserting  a non-overlined part ${\eta}$, and so $t=1$. From the proof above, we see that $\eta\leq g(\pi)<\overline{2\eta}$.

{\noindent Condition (D).}  For any $(k-1)$-band $\{\pi_{i+l}\}_{0\leq l\leq k-2}$ of $\pi$, we wish to show that  $\{\pi_{i+l}\}_{0\leq l\leq k-2}$ is even. There are two cases.

{\noindent Case (D)-1.} $\pi_i>(t+1)\eta$, or $\overline{(t+1)\eta}\leq\pi_i\leq(t+1)\eta$ and $\pi$ is obtain from $\mu$ by inserting   an overlined part $\overline{t\eta}$. Using the condition (A), we know that $\pi$ satisfies the condition (3) in Definition \ref{defi-O-B}. So, we find that $\pi_{i+k-2}\geq t\eta$ with strict inequality if $\pi_i>(t+1)\eta$.

 By construction of $\pi$, we see that $\mu_{i+l}=\pi_{i+l}$ for $0\leq l\leq k-2$. Moreover, $\{\mu_{i+l}\}_{0\leq l\leq k-2}$ is a $(k-1)$-band of $\mu$,  $\overline{V}_\mu(\mu_{i})=\overline{V}_\pi(\pi_{i})$ and $\overline{O}_\mu(\mu_{i+k-2})=\overline{O}_\pi(\pi_{i+k-2})$. Under the condition that $\mu\in\overline{\mathcal{B}}^{\,>}_0(\alpha_1,\ldots,\alpha_\lambda;\eta,k,r)$, we see that the $(k-1)$-band $\{\mu_{i+l}\}_{0\leq l\leq k-2}$ of $\mu$ is even, and so the $(k-1)$-band   $\{\pi_{i+l}\}_{0\leq l\leq k-2}$  of $\pi$ is also even.

{\noindent Case (D)-2.} $\pi_i<\overline{(t+1)\eta}$.  It follows from the condition (B) that $g(\pi)\geq \overline{t\eta}$. By the definition of $g(\pi)$, we have $\pi_i\geq g(\pi)\geq \overline{t\eta}$, and so $\pi_{i+k-2}\geq (t-1)\eta$. It yields that  $\{\pi_{i+l}\}_{0\leq l\leq k-2}$ is a $(k-1)$-band of $\pi$ belonging to $[(t-1)\eta,\overline{(t+1)\eta})$. By Lemmas \ref{t+1} and \ref{t+1-o}, we derive that $\{\mu_{i+l}\}_{0\leq l\leq k-3}$ is a $(k-2)$-band of $\mu$ belonging to $[(t-1)\eta,\overline{(t+1)\eta})$.

Define $\delta(t)=0$ if $\pi$ is obtain from $\mu$ by inserting a non-overlined part ${t\eta}$, or otherwise, $\delta(t)=1$ if $\pi$ is obtain from $\mu$ by inserting an overlined part $\overline{t\eta}$. By the definition of the $(k-1)$-augmentation $C_t$, we have
\begin{equation}\label{mu-equiv-1}
[|\mu_{i}|/\eta]+\cdots+[|\mu_{i+k-3}|/\eta]\equiv \delta(t)+t+r-1+\overline{V}_{\mu}(\mu_{i})+\overline{O}_{\mu}(\mu_{i+k-3})\pmod2.
\end{equation}

Utilizing \eqref{cong-m-band} and \eqref{cong-m-band-o}, we get
\begin{equation*}\label{cong-m-band-delta}
    \begin{split}
&\quad [|\pi_{i}|/\eta]+\cdots+[|\pi_{i+k-2}|/\eta]+\overline{V}_{\pi}(\pi_{i})+\overline{O}_{\pi}(\pi_{i+k-2})\\
&\equiv [|\mu_{i}|/\eta]+\cdots+[|\mu_{i+k-3}|/\eta]+\overline{V}_{\mu}(\mu_{i})+\overline{O}_{\mu}(\mu_{i+k-3})+t+\delta(t)\pmod2.
\end{split}
\end{equation*}
Combining with \eqref{mu-equiv-1}, we arrive at
\begin{equation*}\label{k-2-band-congruence-0000}
[|\pi_{i}|/\eta]+\cdots+[|\pi_{i+k-2}|/\eta]\equiv r-1+\overline{V}_{\pi}(\pi_{i})+\overline{O}_{\pi}(\pi_{i+k-2})\pmod2.
\end{equation*}
This proves that the $(k-1)$-band $\{\pi_{i+l}\}_{0\leq l\leq k-2}$ of $\pi$ is even.

{\noindent Case (D)-3.} $\overline{(t+1)\eta}\leq\pi_i\leq(t+1)\eta$ and $\pi$ is obtain from $\mu$ by inserting   a overlined part ${t\eta}$. With the similar argument in  Case (D)-2, we deduce that $\{\pi_{i+l}\}_{0\leq l\leq k-2}$ is a $(k-1)$-band of $\pi$ belonging to $[(t-1)\eta,{(t+1)\eta}]$.

By the definition of the $(k-1)$-augmentation $C_t$, we see that there exists a $(k-2)$-band   $\{\mu_{s+l}\}_{0\leq l\leq k-3}$ of $\mu$ belonging to $[(t-1)\eta,\overline{(t+1)\eta})$ which is of type N. Using Lemma 3.1, we obtain that $\{\pi_{s+l}\}_{0\leq l\leq k-2}$ is a $(k-1)$-band of $\pi$ belonging to $[(t-1)\eta,\overline{(t+1)\eta})$. With the similar argument in  Case (D)-2, we deduce that the $(k-1)$-band $\{\pi_{s+l}\}_{0\leq l\leq k-2}$ of $\pi$ is even.
Appealing to Lemma \ref{cong-t-0}, we see that the $(k-1)$-band $\{\pi_{i+l}\}_{0\leq l\leq k-2}$ of $\pi$ is even.

By now, we have shown that   $\pi$ is an overpartition
 in $\overline{\mathcal{B}}^{\,=}_0(\alpha_1,\ldots,\alpha_\lambda;\eta,k,r|t)$.  Clearly,  $|\pi|=|\mu|+t\eta$ and $\ell(\pi)=\ell(\mu)+1$. This completes the proof.  \qed

We proceed to show that the $(k-1)$-reduction operation  and  the $(k-1)$-augmentation operation are bijections between $\overline{\mathcal{B}}^{\,>}_0(\alpha_1,\ldots,\alpha_\lambda;\eta,k, r|t)$ and $\overline{\mathcal{B}}^{\,=}_0(\alpha_1,\ldots,\alpha_\lambda;\eta,k, r|t)$.

\begin{thm}\label{theorem4.4}
The $(k-1)$-reduction operation  and  the $(k-1)$-augmentation operation are the inverse map of each other.
\end{thm}

\pf  Let $\pi$ be an overpartition in $\overline{\mathcal{B}}^{\,=}_0(\alpha_1,\ldots,\alpha_\lambda;\eta,k, r|t)$ and let $\mu=D_t(\pi)$. By Lemma \ref{division}, we have $\mu\in\overline{\mathcal{B}}^{\,>}_0(\alpha_1,\ldots,\alpha_\lambda;\eta,k, r|t)$. To prove that $C_t(\mu)=\pi$, we consider the following two cases.

Case 1: $g(\pi)\geq \overline{(t+1)\eta}$. In this case,  we have $s(\pi)=\overline{t\eta}$. Then, $\mu$ is obtained from $\pi$ by removing the overlined part $\overline{t\eta}$. We proceed to show that there are no $(k-2)$-bands of $\mu$ in $[(t-1)\eta,\overline{(t+1)\eta})$. Suppose to the contrary that there exists a  $(k-2)$-band $\{\mu_{s+l}\}_{0\leq l\leq k-3}$ of $\mu$ belonging to $[(t-1)\eta,\overline{(t+1)\eta})$. It follows from the condition (2) in Lemma \ref{t+1-o} that
$\{\pi_{s+l}\}_{0\leq l\leq k-2}$ is a $(k-1)$-band of $\pi$ belonging to $[(t-1)\eta,\overline{(t+1)\eta})$, which contradicts the assumption that $g(\pi)\geq \overline{(t+1)\eta}$. Hence, there are no $(k-2)$-bands of $\mu$ in $[(t-1)\eta,\overline{(t+1)\eta})$. Then, $C_t(\mu)$ is obtained from $\mu$ by inserting an overlined part $\overline{t\eta}$, and so $C_t(\mu)=\pi$.

Case 2: $\overline{t\eta}\leq g(\pi)<\overline{(t+1)\eta}$. Assume that $g(\pi)$ is the $m$-th part of $\pi$. With the similar argument in  Case (D)-2 in the proof of Lemma \ref{combinatio}, we deduce that $\{\pi_{m+l}\}_{0\leq l\leq k-2}$ is a $(k-1)$-band of $\pi$ belonging to $[(t-1)\eta,\overline{(t+1)\eta})$. Under the assumption that $\pi\in\overline{\mathcal{B}}^{\,=}_0(\alpha_1,\ldots,\alpha_\lambda;\eta,k, r|t)$, we have
\begin{equation}\label{k-2-band-congruence-0000-o}
[|\pi_{m}|/\eta]+\cdots+[|\pi_{m+k-2}|/\eta]\equiv r-1+\overline{V}_{\pi}(\pi_{m})+\overline{O}_{\pi}(\pi_{m+k-2})\pmod2.
\end{equation}

By Lemmas \ref{t+1} and \ref{t+1-o}, we obtain that $\{\mu_{m+l}\}_{0\leq l\leq k-3}$ is a $(k-2)$-band of $\mu$ belonging to $[(t-1)\eta,\overline{(t+1)\eta})$.
Define $\delta(t)=0$ if $\mu$ is obtain from $\pi$ by removing a non-overlined part ${t\eta}$, or otherwise, $\delta(t)=1$ if $\mu$ is obtain from $\pi$ by removing an overlined part $\overline{t\eta}$. Using \eqref{cong-m-band} and \eqref{cong-m-band-o}, we get
\begin{equation*}\label{cong-m-band-delta-ooooooooooo}
    \begin{split}
&\quad [|\pi_{m}|/\eta]+\cdots+[|\pi_{m+k-2}|/\eta]+\overline{V}_{\pi}(\pi_{m})+\overline{O}_{\pi}(\pi_{m+k-2})\\
&\equiv [|\mu_{m}|/\eta]+\cdots+[|\mu_{m+k-3}|/\eta]+\overline{V}_{\mu}(\mu_{m})+\overline{O}_{\mu}(\mu_{m+k-3})+t+\delta(t)\pmod2.
\end{split}
\end{equation*}
Combining with \eqref{k-2-band-congruence-0000-o}, we arrive at
\[
[|\mu_{m}|/\eta]+\cdots+[|\mu_{m+k-3}|/\eta]\equiv \delta(t)+t+r-1+\overline{V}_{\mu}(\mu_{m})+\overline{O}_{\mu}(\mu_{m+k-3})\pmod2.
\]

It implies that the $(k-2)$-band $\{\mu_{m+l}\}_{0\leq l\leq k-3}$ of $\mu$ is of type N if $\mu$ is obtain from $\pi$ by removing a non-overlined part ${t\eta}$, or otherwise, $\{\mu_{m+l}\}_{0\leq l\leq k-3}$  is of type O if $\mu$ is obtain from $\pi$ by removing an overlined part $\overline{t\eta}$. By the definition of the $(k-1)$-augmentation $C_t$, we get $C_t(\mu)=\pi$.

In either case, we have shown that $C_t(D_t(\pi))=\pi$. Conversely, let $\mu$ be an overpartition in $\overline{\mathcal{B}}^{\,>}_0(\alpha_1,\ldots,\alpha_\lambda;\eta,k, r|t)$. Using Lemma \ref{combinatio}, we have $C_t(\mu)\in \overline{\mathcal{B}}^{\,=}_0(\alpha_1,\ldots,\alpha_\lambda;\eta,k, r|t)$. It follows from Definitions \ref{defi-division} and \ref{insertion} that $D_t(C_t(\mu))=\mu$. This completes the proof. \qed

\subsection{Proof of Theorem \ref{eqv-main-0}}

In this subsection, we   demonstrate that   Theorem \ref{eqv-main-0} can be
justified by  repeatedly  using   the $(k-1)$-reduction  and   the $(k-1)$-augmentation operations.

{\it \noindent Proof of Theorem \ref{eqv-main-0}.} Let $\pi$ be an overpartition in $\mathcal{\overline{B}}_0(\alpha_1,\ldots,\alpha_\lambda;\eta,k,r)$. We  wish to construct a pair of overpartitions  $\Phi(\pi)=(\tau,\mu)$ in  $\mathcal{D}_\eta \times {\mathcal{B}}_{1}(\alpha_1,\ldots,\alpha_\lambda;\eta,k-1,r)$ such that $|\pi|=|\tau|+|\mu|$ and $\ell(\pi)=\ell(\tau)+\ell(\mu)$.  We consider the following two cases:

 Case 1: There are no $(k-1)$-bands of $\pi$ and  there are no overlined parts divisible by $\eta$ in $\pi$. Then set  $\tau=\emptyset$ and   $\mu=\pi$.  By definition,  we see that  $\mu$ is an overpartition in ${\mathcal{B}}_{1}(\alpha_1,\ldots,\alpha_\lambda;\eta,k-1,r)$. Moreover, $|\pi|=|\tau|+|\mu|$ and $\ell(\pi)=\ell(\tau)+\ell(\mu)$.

 Case 2: There exists a $(k-1)$-band of $\pi$ or an overlined part  divisible by $\eta$ in $\pi$.  Set $b=0$, $\tau^{(0)}=\emptyset$, $\pi^{(0)} =\pi$, and execute  the following procedure. Denote the intermediate pairs by $(\tau^{(0)},\pi^{(0)}),(\tau^{(1)},\pi^{(1)})$, and so on.
  \begin{itemize}

\item[(A)]  Set
     \[t_{b+1}=\min\left\{\left[ |s(\pi^{(b)})|/\eta\right],\left[|g(\pi^{(b)})|/\eta\right]\right\}.\]
Since $s(\pi^{(b)})\geq\overline{\eta}$ or $g(\pi^{(b)})\geq \overline{\eta}$, we find that $t_{b+1}\geq 1$ and
\[\pi^{(b)}\in \overline{\mathcal{B}}^{\,=}_0(\alpha_1,\ldots,\alpha_\lambda;\eta,k,r|t_{b+1}).\]

Applying the $(k-1)$-reduction $D_{t_{b+1}}$ to $\pi^{(b)}$, we get
 \[\pi^{(b+1)}=D_{t_{b+1}}(\pi^{(b)}).\]
Using Lemma \ref{division}, we deduce that
$\pi^{(b+1)}\in \overline{\mathcal{B}}^{\,>}_0(\alpha_1,\ldots,\alpha_\lambda;\eta,k,r|t_{b+1}),$
 \begin{equation}\label{theorem4.1weig}
  |\pi^{(b+1)}|=|\pi^{(b)}|-\eta t_{b+1},
 \end{equation}
 and
  \begin{equation}\label{theorem4.1len}
     \ell(\pi^{(b+1)})=\ell(\pi^{(b)})-1.
 \end{equation}
Then, insert $\eta t_{b+1}$  into $\tau^{(b)}$ as a part to get $\tau^{(b+1)}$.

\item[(B)] Replace $b$ by $b+1$. If there are no  $(k-1)$-bands of $\pi^{(b)}$ and there are no  overlined parts divisible by $\eta$ in $\pi^{(b)}$, then we are done. Otherwise, go back to (A).

 \end{itemize}

 Clearly, the above procedure terminates after at most $\ell(\pi)$ iterations. Assume that it terminates with  $b=m$, that is, there are no  $(k-1)$-bands of $\pi^{(m)}$ and there are no  overlined parts divisible by $\eta$ in $\pi^{(m)}$.  Set
    \[\tau=\tau^{(m)}=(\eta t_m,\ldots,\eta t_1) \quad \text{and} \quad  \mu=\pi^{(m)}.\]

  Then, we have $\mu=\pi^{(m)} \in {\mathcal{B}}_{1}(\alpha_1,\ldots,\alpha_\lambda;\eta,k-1,r)$.
 Using Proposition \ref{daxiao-0-4}, we obtain that for $0\leq b\leq m-2$,
 $t_{b+2}>t_{b+1}\geq 1$. It implies that $ \tau \in \mathcal{D}_\eta$.   Moreover, it is clear from \eqref{theorem4.1weig} and \eqref{theorem4.1len} that  $|\pi|=|\tau|+|\mu|$ and $\ell(\pi)=\ell(\tau)+\ell(\mu)$.
Hence,  $\Phi$ is the desired map from  $\mathcal{\overline{B}}_0(\alpha_1,\ldots,\alpha_\lambda;\eta,k,r)$ to $\mathcal{D}_\eta\times{\mathcal{B}}_{1}(\alpha_1,\ldots,\alpha_\lambda;\eta,k-1,r)$.

To prove that $\Phi$ is a bijection, we shall give the inverse map $\Psi$ of $\Phi$. Let $\tau$ be a partition in $\mathcal{D}_\eta$ and let  $\mu$ be an overpartition in ${\mathcal{B}}_{1}(\alpha_1,\ldots,\alpha_\lambda;\eta,k-1,r)$.
We shall construct  an overpartition $\pi=\Psi(\tau,\mu)\in \mathcal{\overline{B}}_0(\alpha_1,\ldots,\alpha_\lambda;\eta,k,r)$  such that  $|\pi|=|\tau|+|\mu|$ and $\ell(\pi)=\ell(\tau)+\ell(\mu)$.  There are  two cases.

 Case 1:  $\tau=\emptyset$. Then set $\pi=\mu$. It follows from  $\mu\in{\mathcal{B}}_{1}(\alpha_1,\ldots,\alpha_\lambda;\eta,k-1,r)$ that there are no  $(k-1)$-bands  of $\mu$, and so $\pi=\mu\in\mathcal{\overline{B}}_0(\alpha_1,\ldots,\alpha_\lambda;\eta,k,r)$. Moreover, $|\pi|=|\tau|+|\mu|$ and $\ell(\pi)=\ell(\tau)+\ell(\mu)$.

 Case 2: $\tau\neq\emptyset$. Assume that $\tau=(\eta t_m,\eta t_{m-1},\ldots,\eta t_1)$, where $ t_m> t_{m-1}>\cdots> t_1\geq 1$. Starting with $\mu$, apply the $(k-1)$-augmentation repeatedly   to get $\pi$.
  Denote the intermediate overpartitions by $\mu^{(m)},\ldots,\mu^{(0)}$ with $\mu^{(m)}=\mu$ and  $\mu^{(0)}=\pi$. Under the assumption that  $\mu\in{\mathcal{B}}_{1}(\alpha_1,\ldots,\alpha_\lambda;\eta,k-1,r)$, we deduce that there are no  $(k-1)$-bands of $\mu$ and  there are no   overlined parts divisible by $\eta$ in $\mu$. So, we have $s(\mu)=+\infty$ and $g(\mu)=+\infty$, which yields
 $\mu^{(m)}=\mu\in\overline{\mathcal{B}}^{\,>}_0
    (\alpha_1,\ldots,\alpha_\lambda;\eta,k,r| t_{m}).$
 Set  $b=m$, and execute  the following procedure.

\begin{itemize}

\item[(A)] Since
\[\mu^{(b)} \in \overline{\mathcal{B}}^{\, >}_0    (\alpha_1,\ldots,\alpha_\lambda;\eta,k,r|t_{b}),\]
then we apply the the $(k-1)$-augmentation $C_{t_{b}}$ to $\mu^{(b)}$ to get $\mu^{(b-1)}$, that is,
    \[\mu^{(b-1)}=C_{t_{b}}(\mu^{(b)}).\]

Appealing to Lemma \ref{combinatio}, we see that $\mu^{(b-1)}\in \overline{\mathcal{B}}^{\,=}_0(\alpha_1,\ldots,\alpha_\lambda;\eta,k,r|t_{b}),$

 \begin{equation}\label{insertaa-1}
 |\mu^{(b-1)}|=|\mu^{(b)}|+\eta  t_{b},
 \end{equation}
 and
 \begin{equation}\label{insertaa-1t}
\ell(\mu^{(b-1)})=\ell(\mu^{(b)})+1. \end{equation}

   \item[(B)] Replace $b$ by $b-1$. If $b=0$, then we are done. Otherwise,  since $ t_{b+1}> t_{b}$, it follows from Proposition \ref{daxiao-0-4}  that
  \[\mu^{(b)}\in \overline{\mathcal{B}}^{\,>}_0(\alpha_1,\ldots,\alpha_\lambda;\eta,k,r|
   t_{b}).\]
  Go back to (A).
     \end{itemize}

    Eventually, the above process yields an overpartition $\pi=\mu^{(0)}\in  \overline{\mathcal{B}}^{\,=}_0(\alpha_1,\ldots,\alpha_\lambda;\eta,k,r|t_{1})$,
    and so   $\pi$ is an overpartition in $\mathcal{\overline{B}}_0(\alpha_1,\ldots,\alpha_\lambda;\eta,k,r)$.  In light of \eqref{insertaa-1} and  \eqref{insertaa-1t}, we get
\begin{equation*}
     |\pi|=|\mu^{(0)}|=|\mu^{(m)}|+\eta t_{m}+\cdots+\eta t_{1}=|\mu|+|\tau|,
    \end{equation*}
    and
 \begin{equation*}
     \ell(\pi)=\ell(\mu^{(0)})=\ell(\mu^{(m)})+m=\ell(\mu)+\ell(\tau).
    \end{equation*}
Therefore,  $\Psi$ is a map from  $\mathcal{D}_\eta\times{\mathcal{B}}_{1}(\alpha_1,\ldots,\alpha_\lambda;\eta,k-1,r)$ to $\mathcal{\overline{B}}_0(\alpha_1,\ldots,\alpha_\lambda;\eta,k,r)$.
By  Theorem \ref{theorem4.4},  we obtain that $\Phi$ and $\Psi$ are the inverse map of each other. Hence $\Phi$ is a bijection   between $\mathcal{\overline{B}}_0(\alpha_1,\ldots,\alpha_\lambda;\eta,k,r)$ and $\mathcal{D}_\eta\times{\mathcal{B}}_{1}(\alpha_1,\ldots,\alpha_\lambda;\eta,k-1,r)$. This completes the proof.   \qed

 \subsection{An example}

 We conclude this section with an example for  the bijection $\Phi$ in Theorem \ref{eqv-main-0}.  Let \[\pi=(60,60,\overline{53},\overline{50},\overline{47},40,\overline{37},\overline{33},30,\overline{27},\overline{23},20,\overline{20},\overline{10},\overline{7},\overline{3})\]
 be an overpartition in $\overline{\mathcal{B}}_0(3,7;10,5,3)$.
The pair of overpartitions  $\Phi(\pi)=(\tau,\mu)$ is obtained by successively applying the $4$-reduction to $\pi$. The detailed process is given below.

 \begin{itemize}
     \item   Set  $\tau^{(0)}=\emptyset$ and $\pi^{(0)}=\pi$. It can be checked that  $s(\pi^{(0)})=\overline{10}$ and $g(\pi^{(0)})=\overline{27}$. Let
  \[t_{1}=\min\left\{\left[ |s(\pi^{(0)})|/\eta\right],\left[|g(\pi^{(0)})|/\eta\right]\right\}=1.\]
Then, we have $\pi^{(0)}\in\overline{\mathcal{B}}^{\,=}_0(3,7;10,5,3|1)$.  Removing $\overline{10}$ from $\pi^{(0)}$ to get
 \[\pi^{(1)}=D_1(\pi^{(0)})=(60,60,\overline{53},\overline{50},\overline{47},40,\overline{37},\overline{33},30,\overline{27},\overline{23},20,\overline{20},\overline{7},\overline{3}).\]

 Setting $\tau^{(1)}=(10)$ and using Lemma \ref{division}, we obtain that $\pi^{(1)}\in\overline{\mathcal{B}}^{\,>}_0(3,7;10,5,3|1)$.

   \item Since $s(\pi^{(1)})=\overline{20}$ and $g(\pi^{(1)})=\overline{27}$, we have
\[t_2=\min\left\{\left[ |s(\pi^{(1)})|/\eta\right],\left[|g(\pi^{(1)})|/\eta\right]\right\}=2,\]
whence  $\pi^{(1)}\in\overline{\mathcal{B}}^{\,=}_0(3,7;10,5,3|2)$. Removing $\overline{20}$ from $\pi^{(1)}$, we  get
 \[\pi^{(2)}=D_2(\pi^{(1)})=(60,60,\overline{53},\overline{50},\overline{47},40,\overline{37},\overline{33},30,\overline{27},\overline{23},20,\overline{7},\overline{3}).\]
 Setting $\tau^{(2)}=(20,10)$ and using Lemma \ref{division}, we obtain that $\pi^{(2)}\in \overline{\mathcal{B}}^{\,>}_0(3,7;10,5,3|2)$.

\item  Since $s(\pi^{(2)})=\overline{50}$ and $g(\pi^{(2)})=30$,  we have
\[t_3=\min\left\{\left[ |s(\pi^{(2)})|/\eta\right],\left[|g(\pi^{(2)})|/\eta\right]\right\}=3,\]
which implies that $\pi^{(2)}\in\overline{\mathcal{B}}^{\,=}_0(3,7;10,5,3|3)$.  Apply the $4$-reduction to $\pi^{(2)}$ to get $\pi^{(3)}$, namely, $\pi^{(3)}$ is obtained from $\pi^{(2)}$ by  removing a non-overlined part ${30}$. We get
 \[\pi^{(3)}=D_3(\pi^{(2)})=(60,60,\overline{53},\overline{50},\overline{47},40,\overline{37},\overline{33},\overline{27},\overline{23},20,\overline{7},\overline{3}).\]
 Setting $\tau^{(3)}=(30,20,10)$ and using  Lemma \ref{division}, we see that $\pi^{(3)}\in \overline{\mathcal{B}}^{\,>}_0(3,7;10,5,3|3)$.

 \item  Note that $s(\pi^{(3)})=\overline{50}$ and $g(\pi^{(3)})=+\infty$, then we have
 \[t_4=\min\left\{\left[ |s(\pi^{(2)})|/\eta\right],\left[|g(\pi^{(2)})|/\eta\right]\right\}=5,\]
and so $\pi^{(3)}\in\overline{\mathcal{B}}^{\,=}_0(3,7;10,5,3|5)$.
Apply the $4$-reduction to $\pi^{(3)}$ to get $\pi^{(4)}$, namely, $\pi^{(4)}$ is obtained by  removing the overlined part $\overline{50}$ from $\pi^{(3)}$. We get
  \[\pi^{(4)}=D_5(\pi^{(3)})=(60,60,\overline{53},\overline{47},40,\overline{37},\overline{33},\overline{27},\overline{23},20,\overline{7},\overline{3}).\]
 Setting $\tau^{(4)}=(50,30,20,10)$ and using Lemma \ref{division}, we have $\pi^{(4)}\in \overline{\mathcal{B}}^{\,<}_0(3,7;10,5,3|5)$. Eventually,  there are no $4$-bands of $\pi^{(4)}$ and there are no overlined parts divisible by $10$ in $\pi^{(4)}$.
 \end{itemize}
We now get  a  pair of overpartitions $(\tau,\mu)$ with
\begin{equation}\label{omega-reverse}
 \tau=\tau^{(4)}=(50,30,20,10)\text{ and } \mu=\pi^{(4)}=(60,60,\overline{53},\overline{47},40,\overline{37},\overline{33},\overline{27},\overline{23},20,\overline{7},\overline{3})
\end{equation}
such that    $(\tau, \mu) \in \mathcal{D}_{10} \times \overline{\mathcal{B}}_1(3,7;10,4,3)$, $|\pi|=|\tau|+|\mu|$ and $\ell(\pi)=\ell(\tau)+\ell(\mu)$.

Conversely, given   $(\tau, \mu) \in \mathcal{D}_{10} \times \overline{\mathcal{B}}_1(3,7;10,4,3)$ as in \eqref{omega-reverse}, we may recover the overpartition $\pi$ by successively applying the $4$-augmentation
operation. More precisely, the reverse process goes as follows.

\begin{itemize}
    \item Merge $50$ into $\mu^{(4)}=\mu$ to get $\mu^{(3)}$. Since there are no $4$-bands of $\mu^{(4)}$ and there are no overlined parts divisible by $10$ in $\mu^{(4)}$,  we have   $s(\mu^{(4)})=+\infty$ and $g(\mu^{(4)})=+\infty$, which implies that $\mu^{(4)}\in\overline{\mathcal{B}}^{\,>}_0(3,7;10,5,3|5)$. There are no $3$-bands of $\mu^{(4)}$ belonging to $[40,\overline{60})$. Then insert $\overline{50}$ into $\mu^{(4)}$  as an overlined part to  get
     \[\mu^{(3)}=C_{5}(\mu^{(4)})=(60,60,\overline{53},\overline{50},\overline{47},40,\overline{37},\overline{33},\overline{27},\overline{23},20,\overline{7},\overline{3}).\]
     Using Lemma \ref{combinatio}, we obtain that $\mu^{(3)}\in\overline{\mathcal{B}}^{\,=}_0(3,7;10,5,3|5)$.

     \item   Merge $30$ into $\mu^{(3)}$ to get $\mu^{(2)}$. By Proposition \ref{daxiao-0-4}, we have $\mu^{(2)}\in\overline{\mathcal{B}}^{\,>}_0(3,7;10,5,3|3)$. There is a $3$-band $\{\overline{27},\overline{23},20\}$ of $\mu^{(3)}$ belonging to $[20,\overline{40})$ which is of type N. Apply the $4$-augmentation to $\mu^{(3)}$ to get $\mu^{(2)}$, namely, $\mu^{(2)}$ is obtained by  inserting ${30}$  into $\mu^{(3)}$ as a non-overlined part.  We get
     \[\mu^{(2)}=C_{3}(\mu^{(3)})=(60,60,\overline{53},\overline{50},\overline{47},40,\overline{37},\overline{33},{30},\overline{27},\overline{23},20,\overline{7},\overline{3}).\]
In light of Lemma \ref{combinatio}, we deduce that $\mu^{(2)}\in\overline{\mathcal{B}}^{\,=}_0(3,7;10,5,3|3)$.

  \item   Merge $20$ into $\mu^{(2)}$ to get $\mu^{(1)}$.
By Proposition \ref{daxiao-0-4}, we get $\mu^{(2)}\in\overline{\mathcal{B}}^{\,>}_0(3,7;10,5,3|2)$.
Note that $g(\pi^{(2)})=30$, then it follows from the condition (2) in  Proposition \ref{k-1-band-con} that
all $3$-bands of $\mu^{(2)}$ belonging to $[10,\overline{30})$ are of type O. Then  insert $\overline{20}$  into $\mu^{(2)}$ as an overlined part gives
    \[\mu^{(1)}=C_{2}(\mu^{(2)})=(60,60,\overline{53},\overline{50},\overline{47},40,\overline{37},\overline{33},{30},\overline{27},\overline{23},20,\overline{20},\overline{7},\overline{3}).\]
     Using Lemma \ref{combinatio}, we obtain that $\mu^{(1)}\in\overline{\mathcal{B}}^{\,=}_0(3,7;10,5,3|2)$.

    \item Finally, merge $10$ into $\mu^{(1)}$ to get $\mu^{(0)}$. It follows from Proposition \ref{daxiao-0-4} that $\mu^{(1)}\in\overline{\mathcal{B}}^{\,>}_0(3,7;10,5,3|1)$. There are no $3$-bands of $\mu^{(1)}$ belonging to $[0,\overline{20})$. Then insert $\overline{10}$ into $\mu^{(1)}$  as an overlined part to  get
    \[\mu^{(0)}=C_{1}(\mu^{(1)})=(60,60,\overline{53},\overline{50},\overline{47},40,\overline{37},\overline{33},{30},\overline{27},\overline{23},20,\overline{20},\overline{10},\overline{7},\overline{3}).\]
 Using Lemma \ref{combinatio}, we obtain that $\mu^{(0)}\in\overline{\mathcal{B}}^{\,=}_0(3,7;10,4,3|1)$.
\end{itemize}

Set $\pi=\mu^{(0)}$. Then   $\pi$ is an overpartition in $\overline{\mathcal{B}}_0(3,7;10,5,3)$ such that $|\pi|=|\tau|+|\mu|$ and $\ell(\pi)=\ell(\tau)+\ell(\mu)$.

\end{document}